\documentclass[%
11pt,
onecolumn,
tightenlines,
superscriptaddress,
preprintnumbers,
nofootinbib,
 amsmath,amssymb,amsthm,
 aps,
eqsecnum,
]{revtex4-2}

\usepackage{isomath}
\usepackage{amsmath,amsthm}
\usepackage{amsbsy}
\usepackage{amssymb}
\usepackage{amscd}
\usepackage{amsfonts}
\usepackage{stmaryrd}
\usepackage{siunitx}
\usepackage{euscript}
\usepackage[utf8]{inputenc}
\usepackage[T1]{fontenc}
\usepackage{newtxtext} 
\everymath{\displaystyle}
\usepackage{exscale}

\usepackage{graphicx}
\usepackage{boxedminipage}
\usepackage{calc}
\usepackage[usenames,dvipsnames]{xcolor}
\graphicspath{ {media/} }
\usepackage[caption=false]{subfig}

\usepackage{setspace}
\usepackage{enumitem}
\setitemize{noitemsep,topsep=0pt,parsep=0pt,partopsep=0pt}
\setenumerate{noitemsep,topsep=0pt,parsep=0pt,partopsep=0pt}
\setdescription{noitemsep,topsep=0pt,parsep=0pt,partopsep=0pt}

\usepackage{hyperref}

\usepackage[small]{titlesec}

\titlespacing*{\section}{0pt}{12pt plus 4pt minus 2pt}{2pt plus 2pt minus 2pt}
\titlespacing*{\subsection}{0pt}{12pt plus 4pt minus 2pt}{2pt plus 2pt minus 2pt}
\titlespacing*\subsubsection{0pt}{12pt plus 4pt minus 2pt}{2pt plus 2pt minus 2pt}
\titlespacing*\paragraph{0pt}{12pt plus 4pt minus 2pt}{2pt plus 2pt minus 2pt}

\makeatletter

    \renewcommand*{\p@subsection}{}
    
    \renewcommand*{\p@subsubsection}{}
\makeatother

\usepackage{isomath}
\usepackage{amsmath}
\usepackage{amssymb}
\usepackage{amscd}
\usepackage{amsfonts}

\newcommand{\half}{\frac{1}{2}}

\newtheorem{remark}{Remark}[section]

\DeclareMathOperator*{\argmin}{argmin}

\newcommand{\const}{\mathrm{const.}}

\newcommand{\parderiv}[2]{\frac{\partial #1}{\partial #2}}
\newcommand{\dm}{\ \mathrm{d}}

\newcommand{\Wcirc}{\overset{\circ}{W}}

\begin{document}

\preprint{To appear in Journal of the Mechanics and Physics of Solids (\url{https://doi.org/10.1016/j.jmps.2021.104716})}

\title{\Large{Phase-Field Modeling and Peridynamics for Defect Dynamics, \\ and an Augmented Phase-Field Model with Viscous Stresses}}

\author{Janel Chua}
    \email{songlinc@andrew.cmu.edu}
    \affiliation{Department of Civil and Environmental Engineering, Carnegie Mellon University}

\author{Vaibhav Agrawal}
    \altaffiliation{Currently at Intel Corporation.}
    \affiliation{Department of Civil and Environmental Engineering, Carnegie Mellon University}

\author{Timothy Breitzman}
    \affiliation{Air Force Research Laboratory}
    
\author{George Gazonas}
    \affiliation{\mbox{CCDC Army Research Laboratory, Attn: FCDD-RLW-MB, Aberdeen Proving Ground, MD 21005, USA}}

\author{Kaushik Dayal}
    \affiliation{Department of Civil and Environmental Engineering, Carnegie Mellon University}
    \affiliation{Center for Nonlinear Analysis, Department of Mathematical Sciences, Carnegie Mellon University}
    \affiliation{Department of Materials Science and Engineering, Carnegie Mellon University}
    
\date{\today}

%%%%%%%%%%%%%%%%%%%%%
%%%%%%%%%%%%%%%%%%%%%
%%%%%%%%%%%%%%%%%%%%%
%%%%%%%%%%%%%%%%%%%%%
\begin{abstract}
    This work begins by applying peridynamics and phase-field modeling to predict 1-d interface motion with inertia in an elastic solid with a non-monotone stress-strain response.
    In classical nonlinear elasticity, it is known that subsonic interfaces require a kinetic law, in addition to momentum balance, to obtain unique solutions; in contrast, for supersonic interfaces, momentum balance alone is sufficient to provide unique solutions.
    This work finds that peridynamics agrees with this classical result, in that different choices of regularization parameters provide different kinetics for subsonic motion but the same kinetics for supersonic motion.
    In contrast, conventional phase-field models coupled to elastodynamics are unable to model, even qualitatively, the supersonic motion of interfaces.
    This work identifies the shortcomings in the physics of standard phase-field models to be:
    (1) the absence of higher-order stress to balance unphysical stress singularities, and 
    (2) the ability of the model to access unphysical regions of the energy landscape. 

    Based on these observations, this work proposes an augmented phase-field model to introduce the missing physics.
    The augmented model adds:
    (1) a viscous stress to the momentum balance, in addition to the dissipative phase-field evolution, to regularize singularities; and
    (2) an augmented driving force that models the physical mechanism that keeps the system out of unphysical regions of the energy landscape.
    When coupled to elastodynamics, the augmented model correctly describes both subsonic and supersonic interface motion.
    The augmented model has essentially the same computational expense as conventional phase-field models and requires only minor modifications of numerical methods, and is therefore proposed as a replacement to the conventional phase-field models.
\end{abstract}

\maketitle

%%%%%%%%%%%%%%%%%%%%%
%%%%%%%%%%%%%%%%%%%%%
%%%%%%%%%%%%%%%%%%%%%
%%%%%%%%%%%%%%%%%%%%%
\section{Introduction}

Peridynamics \cite{silling2000reformulation} and phase-field modeling \cite{abdollahi_arias,ambati2015review,Chen_Phasefield} are currently the leading approaches to model the evolution of microstructure and defects.
An important open question is whether there are qualitative differences in the predictions of these models that cannot be resolved simply by calibration.
That is, given sufficient calibration of model parameters, can both models provide similar predictions for phenomena of interest?
Or, are there settings in which these models -- irrespective of the sophistication of the calibration -- necessarily provide qualitatively different predictions?
If the predictions are different, which -- if either -- could reasonably be considered to be correct?

We examine this question in the context of 1-d interface motion with inertia in a material with a non-monotone stress-strain response.
Classical elasticity has shown that: (1) a kinetic law is required, in addition to momentum balance, to obtain unique solutions for subsonic motion; (2) in contrast, momentum balance alone is sufficient to provide unique solutions for supersonic motion \cite{abeyaratne2006evolution,truskinovsky1993kinks}.
We find, in brief, that peridynamics agrees with classical elasticity while standard phase-field models do not.
Following \cite{abeyaratne1991kinetic} for strain-gradient models, we find that different choices of regularizing parameters in a given peridynamic model gives rise to different kinetics for subsonic motion but the same kinetics for supersonic motion.
In contrast, we show that standard phase-field models are unable to model, even qualitatively, the supersonic motion of interfaces; supersonic motion is shown to necessarily require unbounded stresses.

Given this clear qualitative difference between peridynamics and phase-field models, the next question is which could be considered more reliable?
While peridynamics and phase-field models can model complex phenomena that are beyond the reach of classical elasticity, it is reasonable to require that these more complex models recover classical elasticity when it is applicable.
We therefore propose an augmentation of phase-field models that agrees with the predictions of classical elasticity.

We highlight that a major area of application of both peridynamics and phase-field modeling is to model microstructure evolution, defect motion, dynamic fracture, and so on.
While these phenomena are far more complex than 1-d interface motion, we focus on the latter for several reasons.
First, there is a clear benchmark solution in 1-d interface motion, unlike more complex phenomena.
Second, if predictions do not agree even in simple settings, they are unlikely to agree in more complex settings.
Third, the reason for the disagreement in predictions is easier to understand in a simple setting.

%%%%%%%%%%%%%%%%%%%%%
%%%%%%%%%%%%%%%%%%%%%
%%%%%%%%%%%%%%%%%%%%%
%%%%%%%%%%%%%%%%%%%%%
\paragraph*{Standard Phase-Field Models Are Unsuitable for Problems with Inertia.}

We next turn to why current phase-field models are largely unsuitable for phenomena in which inertial effects -- rather than energy minimization alone -- play a significant role.

Consider a classical elasticity strain energy density $W(\epsilon)$ that is a nonconvex function of the strain $\epsilon$.
The nonconvexity implies a non-monotone stress-strain response; consequently, there are multiple strain values for a given stress (Figure \ref{fig:stress-strain-response}).
This can lead to the formation of microstructure in which regions of constant strain are separated by singular sharp interfaces.

The sharp interfaces are challenging for numerical calculations, and therefore are typically regularized.
In phase-field models, this regularization is done by introducing a phase-field parameter $\phi$ to keep track of the phase or energy well, and introducing gradients of $\phi$ into the energy to penalize sharp interfaces.
That is, the energy density $W(\epsilon)$ is replaced by $\Wcirc(\epsilon,\phi) = w(\phi) + W_{con} (\epsilon - \epsilon_0(\phi))$, and regularized by adding $|\nabla\phi|^2$.
We notice that while $W(\epsilon)$ is nonconvex in $\epsilon$, $W_{con}(\epsilon,\phi)$ is convex in $\epsilon$ and the nonconvexity is introduced through the nonconvex function $w(\phi)$.

\begin{figure}[htb!]
    \subfloat[]{\includegraphics[width=0.45\textwidth]{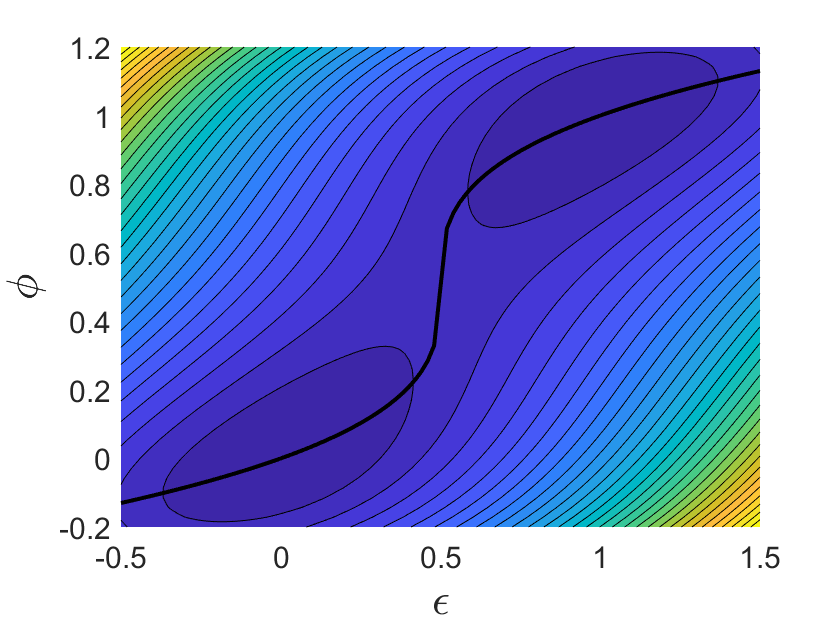}}
    \subfloat[]{\includegraphics[width=0.45\textwidth]{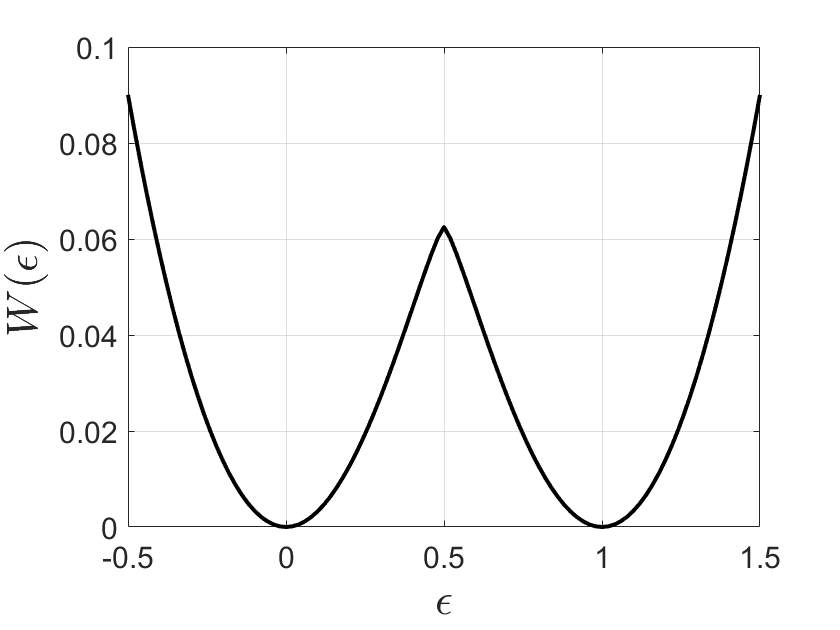}}
    \caption{Left: Energy landscape of a model phase-field energy density $\Wcirc(\epsilon,\phi) =  \phi^2 (1-\phi)^2 + \half (\epsilon-\phi)^2$.
    The bold line shows the value of $\phi$ that minimizes $\Wcirc(\epsilon,\phi)$ for each $\epsilon$, i.e., it plots $\phi_{min}(\epsilon) = \argmin_{\phi} \Wcirc(\epsilon,\phi)$.
    \\
    Right: plot of the energy density along the curve $\phi_{min}(\epsilon)$, i.e. $W(\epsilon):=\Wcirc(\epsilon,\phi_{min}(\epsilon))$, showing multiple minima.}
    \label{fig:phase-field-energy-min}
\end{figure}

The correspondence between classical elasticity and phase-field models can be seen in Figure \ref{fig:phase-field-energy-min}.
We consider homogeneous deformations to enable us to focus on the energy density.
The left panel is a plot of $\Wcirc(\epsilon,\phi)$, and the bold line plots $\phi_{min}(\epsilon)$, where $\phi_{min} = \argmin_{\phi} \Wcirc(\epsilon,\phi)$.
We can then relate the energy densities through $W(\epsilon):=\min_{\phi} \Wcirc(\epsilon,\phi)$.
Appropriate choices of $\Wcirc$ can -- in principle -- be constructed to reproduce a given $W$.

To model the general setting without inertia, it is typical to minimize the total energy with respect to the strain field and use steepest-descent dynamics for the evolution of $\phi$ \cite{clayton2014geometrically,Lun_APL,li2015phase,li2018variational,zhang2016variational,abdollahi_arias,beyerlein2016understanding,vidyasagar2017predicting}.
When inertia is present, it is typical to extremize the action -- or equivalently, to add the inertial term to the momentum balance -- and retain the steepest-descent dynamics for $\phi$ \cite{agrawal2017dependence,bleyer2017dynamic,borden2012phase,ambati2015review,paul2020adaptive,geelen2019phase}.

Energy minimization plays a central role in relating $W$ and $\Wcirc$: it ensures that the material stays near the minimizing curve $\phi_{min}(\epsilon)$ in the energy landscape of Figure \ref{fig:phase-field-energy-min}.
However, in problems that include inertia, energy minimization is not relevant, and the kinetic energy has an important contribution.
The material can then explore parts of the energy landscape far from the minimizing curve; however, the energy landscape away from the minimizing curve has no physical connection to the original classical elasticity energy.
We see in our numerical calculations that the material does explore nonphysical regions of the energy landscape -- particularly at large interface velocities (Section \ref{spfm_kinetics of interface from IVP}).
This is a central reason for phase-field models to fail in correctly modeling supersonic interfaces .
Therefore, an important element of our augmented phase-field model is a physically-motivated driving force that keeps the system away from these nonphysical regions.

%%%%%%%%%%%%%%%%%%%%%
%%%%%%%%%%%%%%%%%%%%%
%%%%%%%%%%%%%%%%%%%%%
%%%%%%%%%%%%%%%%%%%%%
\paragraph*{The Proposed Augmented Phase-Field Model.}

As noted above, one problem with existing phase-field models is that they allow the material to explore unphysical parts of the energy landscape when inertia is important.
A second problem is that there exists a singularity for all supersonic interfaces, i.e., the strain and stress necessarily go to infinity at some point for an interface that moves supersonically.

The first problem shows up in numerical solutions of initial-value problems.
Specifically, as the interface velocity approaches the sonic velocity, we find that the $\phi$ interface has a different velocity and spatial location than the $\epsilon$ interface.
This implies that the spatial region between the two interfaces is in the top-left or bottom-right quadrants of the energy landscape in Figure \ref{fig:phase-field-energy-min} (left) -- the unphysical regions.
The energy $\Wcirc$ in the unphysical regions should be \textit{infinitely} high to respect the original energetic formulation; this would keep the system from exploring these regions, but also cause severe practical difficulties.
However, we notice that if it were possible to set $\Wcirc$ to infinity, this would be reflected in an additional driving force contribution to the dynamical equation for $\phi$.
We therefore augment the dynamical equation for $\phi$ by the missing driving force, and the consequence is that it nudges the system downward in the top-left quadrant, and upward in the bottom-right quadrant.

The second problem is the appearance of an unphysical singularity that forces the stress and strain to go to infinity at a point for a supersonic interface in the standard phase-field formulation (Section \ref{Inability to Model Supersonic Interfaces}).
We also find that a regularization of the momentum equation -- in addition to the usual regularization of phase-field models -- resolves this singularity.
While it possible to use various regularizing stresses, we choose to use a viscous stress because this is the simplest, has a clear physical interpretation, and is readily compatible with standard numerical methods, e.g. FEM with $C^0$ continuity \cite{borden2012phase,ambati2015review, kamensky2018hyperbolic,geelen2019phase}.
We note that all real materials have some level of dissipation, and even in materials in which dissipation is generally small, it can be very important in problems of shocks that are near sonic or supersonic \cite{dafermos2005hyperbolic,abeyaratne2006evolution}.
Therefore, it is not surprising that it plays an important role in phase-field models that aim to be valid when inertial effects are significant.

%%%%%%%%%%%%%%%%%%%%%
%%%%%%%%%%%%%%%%%%%%%
%%%%%%%%%%%%%%%%%%%%%
%%%%%%%%%%%%%%%%%%%%%
\paragraph*{A Note on Strain Gradient Models.}

Strain gradient models are an important class of regularized models of elasticity \cite{abeyaratne1991kinetic,rosakis1995equal,truskinovsky1993kinks,turteltaub1997viscosity}.
They use energetic terms of the form $|\nabla\epsilon|^2$ in the energy to penalize singularly sharp interfaces.
In contrast to phase-field models, they do not introduce any extra fields but have the displacement field as the sole primary field.
However, the strain gradients impose additional restrictions on the continuity of the displacement field that can be challenging for standard FEM.
We notice a heuristic connection between mixed FEM for problems with higher derivatives and the replacement of strain gradient models by phase-field models.
In both cases, we introduce auxiliary variables that nominally relax the smoothness requirements, and then constrain the auxiliary variables to the primary variables.
Phase-field models can be considered analogous to further replacing the constraint by a penalty, which can be justified by energy minimization.

In Section \ref{strain gradient model}, we discuss the findings in the literature on using a strain gradient model to study the problem of interface motion.
In summary, the strain gradient model provides predictions that agree well with classical elasticity: model parameter-dependent kinetics for subsonic motion, and parameter-independent kinetics for supersonic motion.
This raises the question of why one cannot simply use strain gradient models rather than either of peridynamics and phase-field.
While strain gradient models would work well for the particular problem studied here, a number of reasons make strain gradient models unsuitable for broader application.
First, while strain gradient models are useful for regularizing problems that can be described by classical elasticity, it is unclear how to use them for fracture where the displacement itself is discontinuous\footnote{
We mention recent progress by P. Rosakis and coworkers \cite{rosakis2020inverse}.
}.
Second, the dependence of the nucleation and kinetics of interfaces on model parameters is extremely opaque and practically impossible to rationally specify, in contrast to the phase-field models discussed in this paper \cite{agrawal2015dynamic,agrawal2015dynamic-2,alber2005solutions}.
Third, the higher derivatives that appear in the model require nonstandard or restrictive numerical methods, compared to phase-field models that can be solved using standard FEM.

%%%%%%%%%%%%%%%%%%%%%
%%%%%%%%%%%%%%%%%%%%%
%%%%%%%%%%%%%%%%%%%%%
%%%%%%%%%%%%%%%%%%%%%
\paragraph*{Organization.}
Section \ref{sec:formulation} formulates the interface motion problem, and summarizes relevant results from the literature on the solution to the problem in the settings of classical elasticity and strain gradient elasticity.
Sections \ref{sec:peridynamics} and \ref{sec:SPFM} present, respectively, the peridynamic and phase-field model solutions to the interface motion problem.
Section \ref{sec:SPFM-discussion} discusses the reasons that existing phase-field models are unable to model supersonic interface motion.
Section \ref{sec:DPFM} presents and characterizes the augmented phase-field model that is able to correctly model the subsonic and supersonic behavior of interfaces.
Section \ref{sec:discussion} provides further discussion.

%%%%%%%%%%%%%%%%%%%%%
%%%%%%%%%%%%%%%%%%%%%
%%%%%%%%%%%%%%%%%%%%%
%%%%%%%%%%%%%%%%%%%%%
\section{Formulation and Classical Results}
\label{sec:formulation}

The entire paper works in 1-d; our domain is a 1-d bar denoted $\Omega$.
Where closed-form calculations are possible, we will consider $\Omega$ to correspond to the entire real line.
Where numerical calculations are required, we will use a finite bar, and take care to only consider results that are not influenced by the boundaries.

The displacement of the material point at the spatial position $x$ at time $t$ is denoted $u(x,t)$; the stress by $\sigma(x,t)$; the strain by $\epsilon = \partial_x u(x,t)$; and the phase by $\phi(x,t)$.

The material response is formulated to have two phases, denoted phase 1 for the low-strain phase and phase 2 for the high-strain phase.
These phases can coexist in certain situations, and in those situations each phase occupies distinct regions of space separated by interfaces between them.
Throughout the paper, we focus on the motion of individual interfaces.
Various quantities can be discontinuous across the interfaces, but, in line with the fundamental assumptions of continuum mechanics, we require that the displacement is always continuous in space and time.
We denote the location of the interface by $s(t)$ and the velocity by $\dot{s}(t)$.
The jump $g(x=s^+,t) - g(x=s^-,t)$ across the interface of a quantity $g(x,t)$ is denoted $\llbracket g \rrbracket$. 

Our convention is to have phase 1 on the left and phase 2 on the right.
Therefore, $\dot{s} > 0$ corresponds to a transformation of phase 2 to phase 1, and $\dot{s}<0$ corresponds to a transformation of phase 1 to phase 2.

For use in further sections, we introduce $H_l(\cdot)$, a regularized Heaviside / step function that transitions from $0$ to $1$ over a scale $l$ as its argument transitions from negative to positive.
We emphasize that $l$ is not a lengthscale, but is used to scale the phase $\phi$.
For computations, we use the choice $H_l\left(x\right) = \half\left(1+\tanh\left(x/l\right)\right)$.

%%%%%%%%%%%%%%%%%%%%%
%%%%%%%%%%%%%%%%%%%%%
%%%%%%%%%%%%%%%%%%%%%
%%%%%%%%%%%%%%%%%%%%%
\subsection{Material Response}
\label{sec:material-response}

\begin{figure}[htb!]
    \subfloat[]{\includegraphics[width=0.45\textwidth]{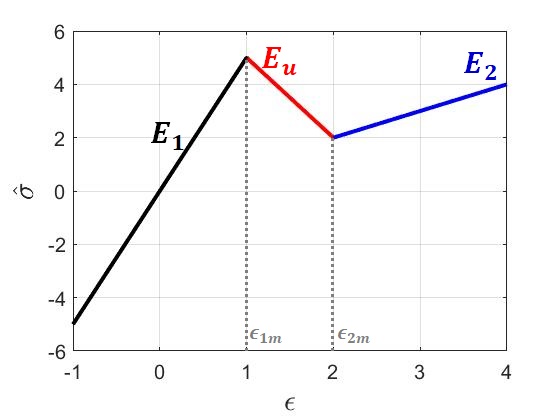}}
    \subfloat[]{\includegraphics[width=0.45\textwidth]{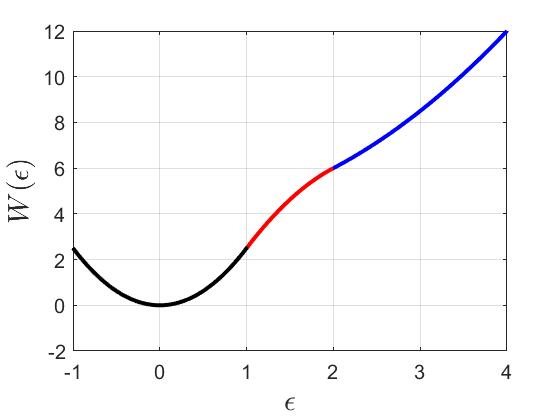}}
    \caption{Material response. (a) Stress-strain curve with $E_1=5$ and $E_2=1$. (b) Strain energy density.}
    \label{fig:stress-strain-response}
\end{figure}

The stress-strain response $\hat\sigma(\epsilon)$ and strain energy density $W(\epsilon)=\int_{0}^{\epsilon} \hat{\sigma}(\tilde{\epsilon}) \dm \tilde\epsilon$ are plotted as a function of strain $\epsilon$ in Figure \ref{fig:stress-strain-response}, and have the following expressions:
\begin{equation}
\label{eqn:stress-energy-response}
\begin{split}
    \hat{\sigma}({\epsilon})
    &= 
    \begin{cases}
        E_1\epsilon           & \text{if } \epsilon \leq \epsilon_{1m}\\
        -E_u\epsilon + C_u    & \text{if } \epsilon_{1m}<\epsilon \leq \epsilon_{2m}\\
        E_2\epsilon           & \text{if } \epsilon>\epsilon_{2m}
    \end{cases}
    \\
    \Rightarrow W({\epsilon}) &= 
    \begin{cases}
        \half E_1\epsilon^2   
        & \text{if } \epsilon \leq \epsilon_{1m}
        \\
        -\half E_u \left(\epsilon^2-\epsilon_{1m}^2\right) + C_u\left(\epsilon-\epsilon_{1m}\right) + \half E_1\epsilon_{1m}^2  
        & \text{if } \epsilon_{1m} <\epsilon \leq \epsilon_{2m}
        \\
        \half E_2\epsilon^2 + 
        \underbrace{
            \left(\half E_1 \epsilon_{1m}^2 -\half E_2 \epsilon_{2m}^2 -\half E_u \left(\epsilon_{2m}^2 - \epsilon_{1m}^2\right) + C_u \left(\epsilon_{2m} - \epsilon_{1m}\right) \right)
        }_{\Delta\Psi}
        & \text{if } \epsilon>\epsilon_{2m}\\
    \end{cases}
\end{split}
\end{equation}
The quantities $E_1$ and $E_2$ are the elastic moduli of phases 1 and 2; $\epsilon_{1m}$ and $\epsilon_{2m}$ are the limits of existence of phases 1 and 2; and $C_u = \frac{\epsilon_{1m}\epsilon_{2m}(E_2-E_1)}{\epsilon_{1m} - \epsilon_{2m}}$ and $E_u=\frac{-(E_1\epsilon_{1m} - E_2\epsilon_{2m})}{\epsilon_{1m} - \epsilon_{2m}}$ are chosen to ensure continuity of $\hat\sigma(\epsilon)$.

We assume for simplicity that the density $\rho$ is constant, and that $E_1 > E_2$.
We define the sonic velocities $c_1 = \sqrt{E_1/\rho}$ and $c_2 = \sqrt{E_2/\rho}$, corresponding to small-amplitude linearized waves in the phases; notice that $c_1 > c_2$.
Subsonic interfaces have $-c_2 < \dot{s} < c_2$, and supersonic interfaces have $c_2 < \dot{s} < c_1$.
Interfaces with $\dot{s} > c_1$ or $ \dot{s} < -c_2$ are not permitted by momentum balance and thermodynamics.

For the numerical calculations, we use $\rho=1,E_2 / \rho=1, E_1 / \rho = 5$.

%%%%%%%%%%%%%%%%%%%%%
%%%%%%%%%%%%%%%%%%%%%
%%%%%%%%%%%%%%%%%%%%%
%%%%%%%%%%%%%%%%%%%%%
\subsection{Results from Classical Elasticity}\label{classical elasticity}

The behavior of interfaces in the material described in Section \ref{sec:material-response} has been studied, among various other topics, in the body of work of Abeyaratne and Knowles \cite{abeyaratne2006evolution}.
We briefly summarize the relevant details of their results here.

First, we consider the quasistatic setting.
The balance of momentum leads to a PDE where the fields are smooth and a jump condition at the interface:
\begin{align}
    & \text{Field equation: } \partial_x \sigma = 0
    \\
    & \text{Jump condition: } \llbracket \sigma \rrbracket = 0
\end{align}
where $\sigma(x,t) = \hat\sigma(\partial_x u)$.

This implies that the stress $\sigma(x,t)$ is constant in the bar, i.e. $\sigma(x,t) = \sigma_0(t)$.
Considering load control, i.e. $\sigma_0(t)$ specified, it is clear from Figure \ref{fig:stress-strain-response} that when $\sigma_0(t) \in [2,5]$, there are an infinity of solutions that satisfy equilibrium.
Specifically, {\em any} displacement field that everywhere has derivative $\partial_x u \in \{\sigma_0/E_1, \sigma_0/E_2\}$ will satisfy momentum balance.
In the context of displacement control, the situation is similar when the average strain in the bar, $\frac{u(L,t)-u(0,t)}{L}$, is in the range $[\frac{2}{5},5]$.

Each discontinuity in $\partial_x u$ corresponds to an interface across which the strain jumps.
Even in the simplified setting of a single interface, the nonuniqueness persists; for instance, given $\sigma_0(t)$, the solution has the form\footnote{We could also have phase 1 on the right and phase 2 on the left.}:
\begin{equation}
\label{eqn:qstatic-soln}
    \partial_x u(x,t) = 
    \begin{cases} 
        \sigma_0/E_1 & \text{ if } x < s(t) \\ 
        \sigma_0/E_2 & \text{ if } x \ge s(t)
    \end{cases}
\end{equation}
where $s(t)$ can be arbitrary.
The nonuniqueness in the quasistatic setting is typically resolved using energy minimization as a selection mechanism \cite{ericksen1975equilibrium}.
For instance, for solutions of the form in \eqref{eqn:qstatic-soln}, we find $s$ by minimizing the potential energy over $s$.

\begin{figure}[htb!]
    {\includegraphics[width=0.45\textwidth]{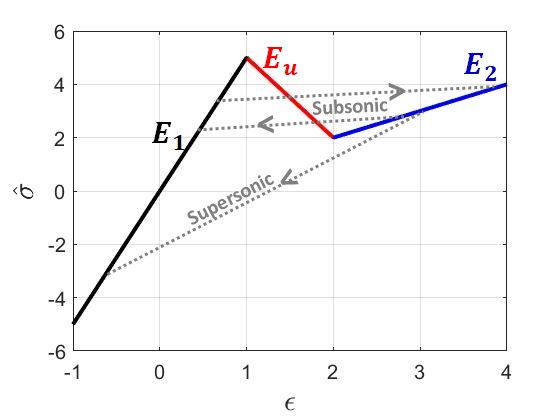}}
    \caption{The stress-strain curve and the ``chord condition''. The chords link the states on either side of the interface.  The slope of a chord is proportional to the square of the interface velocity, and the slope of a stress-strain branch is proportional to the square of the sonic speed for that branch. The subsonic chords have slope such that $|\dot{s}|<c_2$, and the supersonic chord has $c_2 < \dot{s} < c_1$. Notice that it is not possible to construct a chord that has $\dot{s} > c_1$.}
    \label{fig:stress-strain-chord construction}
\end{figure}

Next, we consider the dynamic setting.
The balance of momentum reads:
\begin{align}
    & \text{Field equation: } \partial_x \sigma = \rho \partial_{tt} u
    \label{eqn:dynamic-field}
    \\
    & \text{Jump condition: } \llbracket \sigma \rrbracket = \rho \dot{s}^2 \llbracket \partial_x u \rrbracket 
    \label{eqn:dynamic-jump}
\end{align}
The jump condition \eqref{eqn:dynamic-jump} has an insightful graphical interpretation (Figure \ref{fig:stress-strain-chord construction}).
Writing it as $\dot{s}^2 = \frac{1}{\rho}\frac{\llbracket \sigma \rrbracket}{\llbracket \partial_x u \rrbracket}$, we see that the velocity of the interface is related to the slope of the chord on the stress-strain curve connecting the stress- and strain- states $(\epsilon^-,\sigma^-)$ and $(\epsilon^+,\sigma^+)$ on either side of the interface.
Further, we recall that the sonic velocity in each phase is similarly related to the slope of the stress-strain curve of the corresponding phase, i.e. $c_1^2 = \frac{E_1}{\rho}$ and $c_2^2 = \frac{E_2}{\rho}$.
Indeed, sonic waves are governed by precisely the jump condition \eqref{eqn:dynamic-jump}, but with both end states on the same branch. 

The graphical interpretation tells us that any interface that connects the two branches is represented by a chord whose slope must be less than the slope of the phase 1, implying $\dot{s} < c_1$ always.
However, we notice that $\dot{s}$ can be smaller or larger than $c_2$, and we therefore consider two regimes.
First, the subsonic regime, wherein $1 < M_2 := \frac{\dot{s}}{c_2} \leq 1$ and the interface is subsonic with respect to phase 2; and, second, the supersonic regime, wherein $M_2 > 1$ and the interface is supersonic with respect to phase 2.

The subsonic regime $-1 < M_2\leq 1$ inherits the non-uniqueness that was evident in the quasistatic setting.
In particular, initial-boundary-value problems with subsonic interfaces do not have unique solutions; further, since the problem is dynamic, energy minimization is not an appropriate selection principle.
The non-uniqueness can be related to a lack of information of the kinetics of the interface.
That is, we require a kinetic relation that relates $\dot{s}$ to $(\epsilon^-,\sigma^-)$ and $(\epsilon^+,\sigma^+)$ to obtain unique solutions\footnote{The general problem requires also nucleation criteria, but we focus on the behavior of a single already-nucleated interface throughout this paper.}.

The kinetic relation $\dot{s} = \hat{v}(f)$ relates the velocity of the interface to the \emph{driving force} $f$ acting on the interface, and was introduced by \cite{abeyaratne1990driving,heidug1985thermodynamics,truskinovskii1982equilibrium}.
The driving force is given by the expression:
\begin{equation}
\label{eqn:classical-driving-force}
    f = \llbracket W(\epsilon) \rrbracket - \frac{\hat\sigma(\epsilon^+) + \hat\sigma(\epsilon^-)}{2} \llbracket \epsilon \rrbracket 
\end{equation}
It contains information about the state of the material on both sides of the interface, and is precisely the work conjugate of $\dot{s}$.

The supersonic regime $M_2 > 1$, however, has unique solutions without a kinetic relation \cite{abeyaratne1991implications,truskinovsky1993kinks,trofimov2010shocks}.
That is, \eqref{eqn:dynamic-field} and \eqref{eqn:dynamic-jump} have a unique solution when the interface is supersonic with respect to phase 2.
No additional kinetic relation is required, and using such a kinetic relation will generally over-constrain the problem such that there are no admissible solutions. 

We mention that the chord construction -- i.e., momentum balance -- does not rule out interfaces with $M_2 < -1$, but these would have negative dissipation and hence are ruled out by thermodynamics.

%%%%%%%%%%%%%%%%%%%%%
%%%%%%%%%%%%%%%%%%%%%
%%%%%%%%%%%%%%%%%%%%%
%%%%%%%%%%%%%%%%%%%%%
\subsection{Results from Strain Gradient (Viscosity-Capillarity) Models}\label{strain gradient model}

Strain gradient models, also called viscosity-capillarity models, regularize the sharp interfaces of classical elasticity by adding a strain gradient term to account for the surface energy, and a viscous term to account for the dissipation associated with defect motion \cite{fuaciu2006longitudinal,abeyaratne1991implications,abeyaratne2006evolution,rosakis1995equal,truskinovsky1993kinks,turteltaub1997viscosity}.
We use $\sigma(x,t) = \hat\sigma(\partial_x u) + \rho\nu\partial_{xt}u - \rho\lambda \partial_{xxx} u$ to find the field equation of momentum balance:
\begin{equation}
\label{eqn:strain-gradient}
     \partial_x \hat\sigma(\partial_x u) + \rho\nu\partial_{xxt}u - \rho\lambda \partial_{xxxx} u = \rho \partial_{tt}u
\end{equation}
Here, $\nu$ is the coefficient of dissipation and $\lambda$ is the coefficient of surface energy.
The solution is sufficiently smooth due to the higher derivatives, and therefore the jump condition, \eqref{eqn:dynamic-jump}, is not required.
This model typically has unique solutions given appropriate initial and boundary conditions \cite{abeyaratne1991implications}.

An important finding in \cite{abeyaratne1991implications,truskinovsky1993kinks} is that this regularization preserves the key distinction between the subsonic and supersonic regimes.
That is, the kinetics of subsonic interfaces depends sensitively on the choice of $\nu$ and $\lambda$, whereas the kinetics of supersonic interfaces is relatively insensitive to this choice.
One can therefore think of $\nu$ and $\lambda$ as inducing a kinetic relation when it is required for uniqueness, and providing merely a minor regularizing effect when the kinetic relation is not required for uniqueness.

%%%%%%%%%%%%%%%%%%%%%
%%%%%%%%%%%%%%%%%%%%%
%%%%%%%%%%%%%%%%%%%%%
%%%%%%%%%%%%%%%%%%%%%
\subsubsection{Numerical Computation of Kinetics Using Traveling Waves}\label{strain gradient model_traveling waves}

As mentioned above, the addition of strain gradient and dissipation terms effectively induces a kinetic relation.
While the kinetic relations induced by \eqref{eqn:strain-gradient} have been computed in closed-form in \cite{abeyaratne1991implications}, we nonetheless compute these numerically here as a means to both describe as well as verify our numerical scheme, that we will apply later on to other models studied in this paper.

We begin by assuming a traveling-wave form for the solution: $u(x,t) = U(x-\dot{s}t)$.
Using this in \eqref{eqn:strain-gradient}, we find:
\begin{equation}
\label{eqn:strain-gradient-tw}
     \left(\hat\sigma(U')\right)' 
     - \dot{s}\rho\nu U''' - \rho\lambda U'''' = \dot{s}^2 \rho U''
\end{equation}
where $\tilde{x} := x - \dot{s} t$ is the traveling coordinate, and primes represent differentiation with respect to $\tilde{x}$.

We can immediately integrate \eqref{eqn:strain-gradient-tw} once to get:
\begin{equation}
\label{eqn:strain-gradient-tw-3}
     \frac{1}{E_2} \hat\sigma(\mathcal{E}) - M_2 \frac{\nu}{c_2} \mathcal{E}' - \frac{\lambda}{c_2^2} \mathcal{E}'' 
     = 
     M_2^2 \mathcal{E} + C
\end{equation}
where we have nondimensionalized; introduced the strain $\mathcal{E}(\tilde{x}) := U'(\tilde{x})$; and $C$ is the undetermined constant of integration.

The displacement $U$ is not unique due to rigid-body translations, and we have no boundary conditions to fix this; therefore, we will solve \eqref{eqn:strain-gradient-tw-3} directly for the strain for various given values of $M_2$.

We consider possible boundary conditions for \eqref{eqn:strain-gradient-tw-3}.
To do this, consider the limits $\mathcal{E}(\tilde{x}\to\pm\infty)$, denoting these by $\mathcal{E}^{\pm\infty}$.
Using that the derivatives of $\mathcal{E}$ vanish far from the interface \cite{abeyaratne1991implications}, we find the equations:
\begin{equation}
\label{eqn:strain-gradient-tw-4}
     \frac{1}{E_2} \hat\sigma(\mathcal{E}^{+\infty}) = M_2^2 \mathcal{E}^{+\infty} + C
     , \quad 
     \frac{1}{E_2} \hat\sigma(\mathcal{E}^{-\infty}) = M_2^2 \mathcal{E}^{-\infty} + C
\end{equation}
Subtracting the equations above, we find 
\begin{equation}
\label{eqn:strain-gradient-tw-5}
     \frac{\hat\sigma(\mathcal{E}^{+\infty}) - \hat\sigma(\mathcal{E}^{-\infty})}{\mathcal{E}^{+\infty}-\mathcal{E}^{-\infty}}
     =
     E_2 M_2^2
\end{equation}
This equation is precisely the strain-gradient version of \eqref{eqn:dynamic-jump}, and does not provide any new information in the regularized context.
However, we notice that if we use $\mathcal{E}^{+\infty}$ as given data, along with $M_2$ given, we can use \eqref{eqn:strain-gradient-tw-5} to find $\mathcal{E}^{-\infty}$, and vice-versa.
That is, if we pick either of the far-field strains as data along with given $M_2$, we can solve for the other far-field strain, and consequently solve also for the driving force:
\begin{equation}
\label{eqn:strain-grad-f}
    f = W(\mathcal{E}^{+\infty}) - W(\mathcal{E}^{-\infty}) + \frac{\hat{\sigma}(\mathcal{E}^{+\infty}) +\hat{\sigma}(\mathcal{E}^{-\infty}) }{2}\left(\mathcal{E}^{+\infty}-\mathcal{E}^{-\infty}\right) 
\end{equation}
Thus, if we specify $M_2$ and either of $\mathcal{E}^{\pm\infty}$, we can find the driving force and consequently the kinetic relation corresponding to $M_2$ without using \eqref{eqn:strain-gradient-tw} at all; in particular, the dissipative and surface energy contributions play no role in determining the kinetic relation.
We conclude that specifying either of $\mathcal{E}^{\pm\infty}$ overconstrains the problem at a given $M_2$.

Therefore, we do not specify either of $\mathcal{E}^{\pm\infty}$ and treat both as unknown quantities.
Following \cite{agrawal2015dynamic,agrawal2015dynamic-2,Kaushik-peri}, we solve \eqref{eqn:strain-gradient-tw-3} by treating the entire function $\mathcal{E}(\tilde{x})$ as well as $C$ as unknown and solve for them using a least-squares approach.
Defining the residue:
\begin{equation}
\label{eqn:strain-grad-residue-1}
    r(\tilde{x}) := \frac{1}{E_2} \hat\sigma(\mathcal{E}) - M_2 \frac{\nu}{c_2} \mathcal{E}' - \frac{\lambda}{c_2^2} \mathcal{E}'' - M_2^2 \mathcal{E} - C
\end{equation}
we solve by minimizing $R = \int_{-\infty}^{+\infty} r(\tilde{x})^2 \dm\tilde{x}$ over $C$ and functions $\mathcal{E}$.

The computational discretization uses a finite segment of the bar of length $L$ in the translating frame $\tilde{x}$.
We use finite differences and divide the domain into $N$ equal intervals, each of length $\Delta\tilde{x} = L/N$, with nodes $\tilde{x}_i$, and define $\mathcal{E}_i:=\mathcal{E}(\tilde{x}_i)$.
We approximate the residue by the quadrature:
\begin{equation}
\label{eqn:strain-grad-residue-2}
    R\left(\mathcal{E}_1, \mathcal{E}_2, \ldots, \mathcal{E}_{N+1}, C\right) 
    \approx
    \Delta \tilde{x} \sum_{i=2}^N 
    \left(
        \frac{1}{E_2} \hat\sigma(\mathcal{E}_i) 
        - M_2 \frac{\nu}{c_2} \frac{\mathcal{E}_{i+1} - \mathcal{E}_{i-1}}{2\Delta\tilde{x}} 
        - \frac{\lambda}{c_2^2} \frac{\mathcal{E}_{i+1} - 2\mathcal{E}_{i} + \mathcal{E}_{i-1}}{{\Delta\tilde{x}}^2}
        -M_2^2 \mathcal{E}_i 
        - C
     \right)^2 
\end{equation}
and use left- or right- rather than centered- differences at the edge of the domain.
We minimize $R$ in \eqref{eqn:strain-grad-residue-2} over the finite set of variables $\{\mathcal{E}_i, i=2\ldots N\}$ and $C$ to find the traveling wave profile, using a standard monolithic solver.
To test that we are not stuck at a local minimum, we check that the value of $R$ is close to $0$ after the minimization is complete.

\begin{figure}[htb!]
    {\includegraphics[width=0.6\textwidth ]{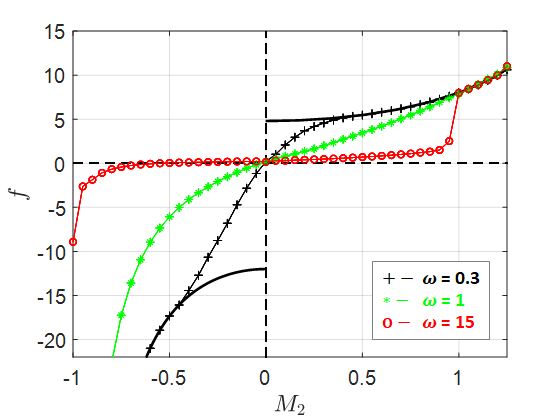}}
    \caption{The kinetic relations in the strain gradient model, for different choices of dissipation and strain gradient coefficients. The kinetic relations vary widely for $-1 < M_2 < 1$, but collapse to a single curve for $M_2 > 1$, in agreement with the predictions of classical elasticity.}
    \label{fig:strain-gradient-kinetics-1}
\end{figure}

\begin{figure}[htb!]
    \subfloat[]{\includegraphics[width=0.33\textwidth]{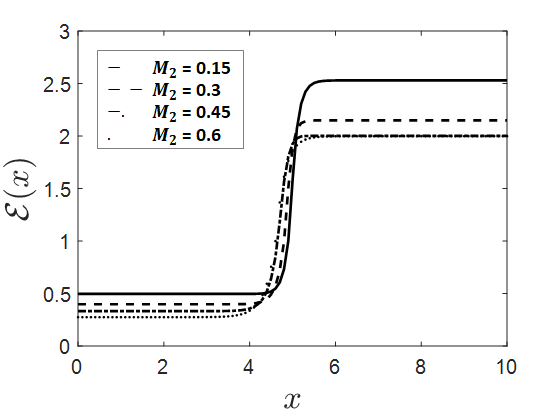}}
    \subfloat[]{\includegraphics[width=0.33\textwidth]{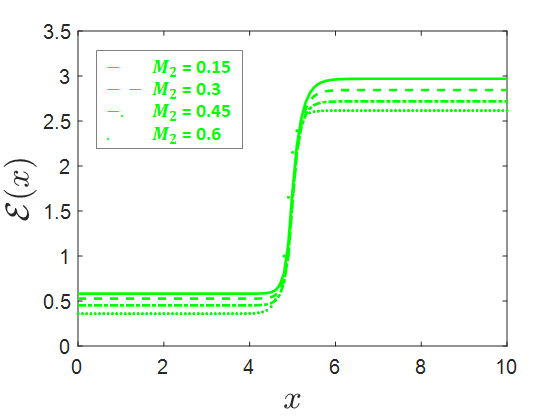}}
    \subfloat[]{\includegraphics[width=0.33\textwidth]{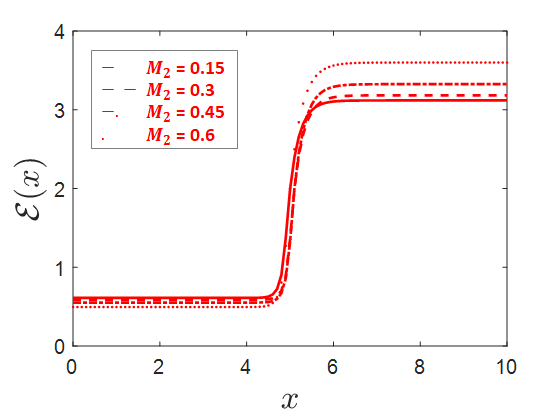}}
    \caption{Representative strain profiles (for different velocities) from the traveling wave computations in the strain gradient model for different choices of dissipation and surface energy coefficient. (a) $\omega = 0.3$ (b) $\omega = 1$ (c) $\omega = 15$.}
    \label{fig:strain-gradient-kinetics-2}
\end{figure}

Using this numerical procedure, for various choices of $M_2$, we compute the strain profile and use this to infer the driving force through \eqref{eqn:strain-grad-f}; this procedure is repeated for several choices of $\omega$ where $\omega = \frac{2\sqrt{\lambda}}{\nu}$.
The corresponding kinetic curves and some representative strain profiles are shown in Figures \ref{fig:strain-gradient-kinetics-1} and \ref{fig:strain-gradient-kinetics-2}.
These match well with the closed-form expressions obtained in \cite{abeyaratne1991implications} and provides us with confidence in the numerical scheme.

%%%%%%%%%%%%%%%%%%%%%
%%%%%%%%%%%%%%%%%%%%%
%%%%%%%%%%%%%%%%%%%%%
%%%%%%%%%%%%%%%%%%%%%
\section{Interface Motion in Peridynamics}
\label{sec:peridynamics}

\newcommand{\bondf}{\mathsf{f}}

Following the model postulated by \cite{silling2000reformulation}, the 1-d peridynamic equation of motion is:
\begin{equation}
\label{eqn:1D peri eqn of motion}
    \rho\partial_{tt}u(x,t) 
    = 
    \int_{\bar x\in\mathbb{R}} \bondf(u(\bar{x},t) - u(x,t),\bar{x}-x) \dm\bar{x}
\end{equation}
where $\bondf(\delta u,\delta x)$ is the bond force between two volume elements with separation in the reference $\delta x:= \bar{x} - x$ and relative displacement $\delta u := u(\bar{x},t) - u(x,t)$.

It is useful to have a dissipative mechanism to account for the dissipation in interface motion.
Instead of adding terms containing the strain rate as in continuum mechanics -- which would nullify the goal of peridynamics of eliminating spatial derivatives -- we follow \cite{Kaushik-peri} and add a dissipative contribution to $\bondf$:
\begin{equation}
    \bondf(\delta u, \delta x)
    = 
    \frac{4}{\sqrt{\pi}} \frac{\delta x}{l_0^3} e^{-\left(\frac{\delta x}{l_0}\right)^2} 
    \left(
        \hat\sigma\left(\frac{\delta u}{\delta x}\right) 
        +
        \frac{\nu_b l_0}{c_2 \delta x} \partial_t\delta u 
    \right)
   \label{eqn:1D peri eqn of motion-f}
\end{equation}
where $l_0$ is the nonlocal length scale and $\nu_b$ is the dimensionless coefficient of viscous bond-level damping.
The argument for $\hat\sigma$ (the classical stress-response function from \eqref{eqn:stress-energy-response}) is the bond strain $\delta u / \delta x$, rather than the classical strain $\epsilon$.
This choice ensures that the stress-strain response for homogeneous deformations, computed sufficiently far from the boundaries, is identical to that chosen in strain-gradient and classical elasticity \cite{weckner2005effect, breitzman2018bond, Kaushik-peri,tupek2014extended}.

We noticed in strain gradient models that changing the parameters $\nu$ and $\lambda$ induced different kinetic relations for subsonic interfaces.
In peridynamics, we will analogously change the parameters $\nu_b$ and $l_0$ to induce different kinetic relations.
Regardless of the values chosen for $\nu_b$ and $l_0$, the form of the expression for the bond force in \eqref{eqn:1D peri eqn of motion-f} gives us the stress-strain response $\hat\sigma$ in the setting of homogeneous deformations.
That is, changing $\nu_b$ and $l_0$ leaves the homogeneous stress-strain response unchanged.

\subsection{Numerical Computation of Kinetics Using Traveling Waves}

Similar to the strain gradient approach, we seek a solution in the form of a traveling wave
\begin{equation}
\label{eqn:traveling wave_peri}
    u(x,t) = U(x-\dot{s}t) = U(\tilde{x})
\end{equation}
Substituting \eqref{eqn:traveling wave_peri} into \eqref{eqn:1D peri eqn of motion}, we have:
\begin{equation}
\label{eqn:1D peri eqn of motion-2}
    M_2^2 U''(\tilde x) = 
    \frac{1}{E_2}\int_{0}^{L} \bondf(U(\bar{\tilde{x}}) - U(\tilde{x}),\bar{\tilde{x}}-\tilde{x}) \dm\bar{\tilde{x}}
\end{equation}
We highlight that we use a finite domain $[0,L]$ for the numerical calculations, and set $L \gg l_0$ and ensure that the interface is far the boundaries.

The residue is defined as:
\begin{equation}
    r(\tilde{x})
    :=
    M_2^2 U''(\tilde{x}) 
    -
    \frac{1}{E_2}\int_{0}^{L} \bondf(U(\bar{\tilde{x}}) - U(\tilde{x}),\bar{\tilde{x}}-\tilde{x}) \dm\bar{\tilde{x}}
\end{equation}

Here, we diverge from the method used for the strain gradient model due to the different nature of the boundary conditions in peridynamics.
Specifically, boundary conditions in peridynamics are applied over a finite layer \cite{silling2000reformulation}.
In the context of a traveling wave where we do not have boundary conditions to apply, the procedure to deal with the boundaries follows \cite{Kaushik-peri,dayal2017leading}.
First, we decompose our domain $[0,L]$, denoted $\Omega$, into boundary regions on the left and right, denoted $\Omega^-$ and $\Omega^+$ respectively, and an interior region, denoted $\mathcal{I}$.
The boundary regions have size $\gg l_0$ such that there are negligible interactions between points in $\mathcal{I}$ and points beyond the boundary regions.
We now define our solution as the minimization of $R = \int_{\tilde x \in \mathcal{I}} r(\tilde{x})^2\dm\tilde{x}$.

We notice that $R$ is a functional of $U(\tilde x)$ over $\Omega$ and not $\mathcal{I}$, despite the domain of integration in the definition of $R$; due to the nonlocal interactions, $R$ involves a double integration.

We now discretize $\Omega$ into $N$ equal intervals, each of length $\Delta\tilde{x} = L/N$, with nodes $\tilde{x}_i$, and define $U_i = U(\tilde{x}_i)$.
Then, we approximate the integrations using the quadrature:
\begin{align}
    \int_{\mathcal{I}} r(\tilde{x})^2\dm\tilde{x} 
    & \approx 
    \Delta \tilde{x} \sum_{\tilde{x}_j \in \mathcal{I}} r(\tilde{x}_j)^2
    \\
    \begin{split}
        r(\tilde{x}_j)
        &=
        M_2^2 U''(\tilde{x}_j)
        - \frac{1}{E_2}\int_\Omega \bondf(U(\bar{\tilde{x}}) - U(\tilde{x}_j),\bar{\tilde{x}}-\tilde{x}_j) \dm\bar{\tilde{x}} 
        \\
        & \approx
        M_2^2 \frac{U_{j+1} - 2U_{j} + U_{j-1}}{\Delta \tilde{x}^2}
        - \frac{\Delta \tilde{x}}{E_2} \sum_{\tilde{x}_i \in \Omega}
        \bondf(U_i - U_j,\tilde{x}_i-\tilde{x}_j)
    \end{split}
\end{align}
and use standard central differences for the derivatives.
We use $\Delta x = \frac{l_0}{10}$ as numerical experiments show that the results are essentially converged with this level of discretization.
We set the average of $U$ over the domain to be $0$ to fix the rigid translation.

We highlight that we minimize the error only over the interior $\mathcal{I}$, but the variables over which we minimize are \textit{all} the nodal values.
We refer to \cite{agrawal2015dynamic,Kaushik-peri,dayal2017leading} for more discussion of this approach.

\begin{figure}[htb!]
    \includegraphics[width=0.6\textwidth]{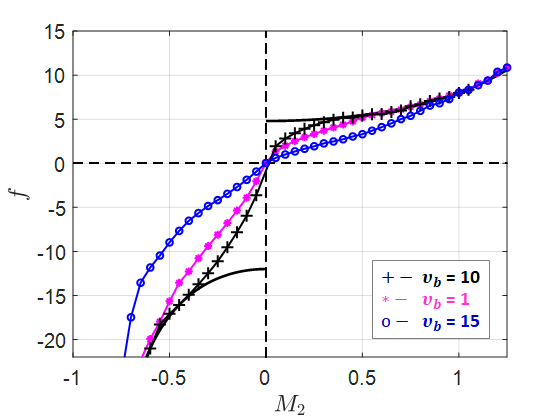}
    \caption{The kinetic relations in peridynamics, for different choices of nonlocality $l_0$ and bond dissipation $\nu_b$. the kinetic relations vary widely for $-1 < M_2 < 1$, but collapse to a single curve for $M_2 > 1$, in agreement with the predictions of classical elasticity.}
    \label{fig:peridynamics-kinetics}
\end{figure}

\begin{figure}[htb!]
    \subfloat[]{\includegraphics[width=0.33\textwidth]{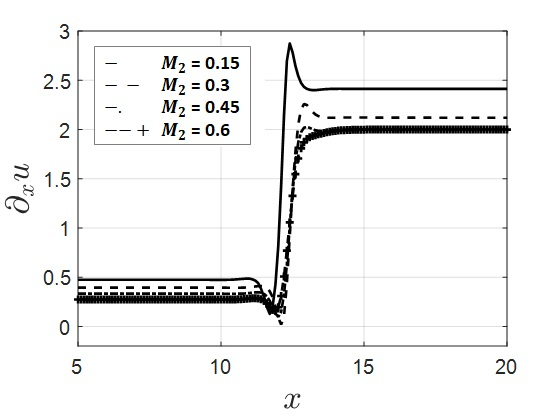}}
    \subfloat[]{\includegraphics[width=0.33\textwidth]{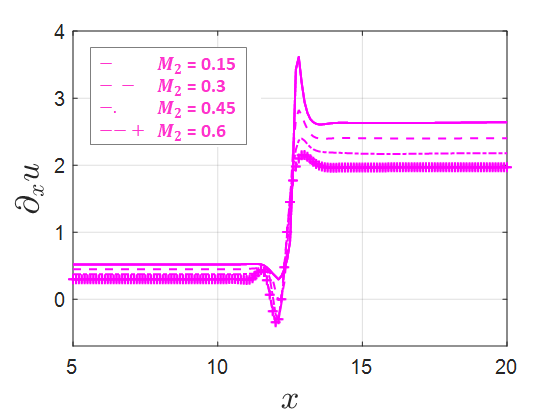}}
    \subfloat[]{\includegraphics[width=0.33\textwidth]{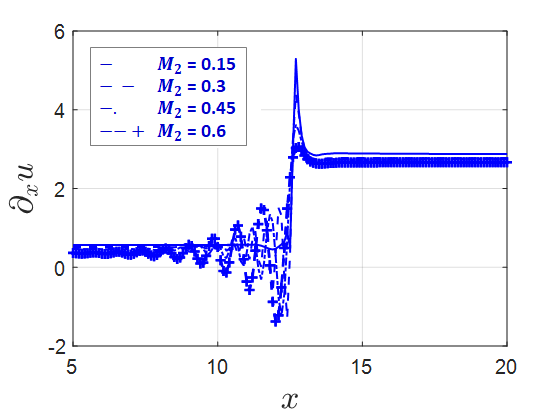}}
    \caption{Representative displacement derivative profiles (for different velocities) from the traveling wave computations in peridynamics for different choices of bond dissipation coefficient. (a) $\nu_b = 10$, (b) $\nu_b = 5$ , (c) $\nu_b = 1$.}
    \label{fig:peridynamics-kinetics-2}
\end{figure}

We use this numerical procedure to compute the traveling wave profile for various values of $M_2$.
Given the displacement profile, we can compute the driving force using \eqref{eqn:strain-grad-f}, thus giving a kinetic relation that relates the interface velocity to the driving force.
This procedure is repeated for several choices of $\nu_b$ and $l_0$.
The corresponding kinetic curves and some representative displacement derivative profiles are shown in Figures \ref{fig:peridynamics-kinetics} and \ref{fig:peridynamics-kinetics-2}.

%%%%%%%%%%%%%%%%%%%%%
%%%%%%%%%%%%%%%%%%%%%
%%%%%%%%%%%%%%%%%%%%%
%%%%%%%%%%%%%%%%%%%%%
\section{Interface Motion in Existing Phase-Field Models}
\label{sec:SPFM}

Two existing phase-field models will be studied in this paper.
While there are several differences between these models, they share 2 key features:
(1) both have gradient regularizations of the phase-field parameter, and not in the momentum balance;
(2) neither has a mechanism to prevent the system from accessing unphysical regions of the energy landscape.
These features are common to all phase-field formulations that we are aware of.

The first model -- or closely-related variations -- is completely standard, and has been used in the overwhelming majority of prior works, e.g. \cite{Chen_Phasefield,Lun_APL,abdollahi_arias,agrawal2017dependence,beyerlein2016understanding,peng20203d,li2015phase,li2018variational,zhang2016variational,albrecht2020phase,collet2020variational}; we will refer to it as the ``standard phase-field model'' for short.
In standard phase-field models, the energy landscape is formulated based on equilibrium principles of energy minimization, and has the overall structure of the form shown in Figure \ref{fig:phase-field-energy-min}.
While the model gives rise to unique evolution of microstructure -- i.e., nucleation and kinetics of interfaces is contained in the model -- the relation between the model parameters and the nucleation and kinetics of interfaces is completely opaque.
Further, nucleation and kinetics are also coupled in the sense that changing the model parameters typically changes both the nucleation and kinetic behavior simultaneously.

The second phase-field model that we will study was motivated by the goal of overcoming the shortcomings discussed just above; we will refer to it as the ``dynamic phase-field model'' for short.
Specifically, the dynamic phase-field model is formulated to provide a transparent separation between energetics, kinetics, and nucleation.
There are distinct model parameters that correspond to each of these, i.e., we can specify the equilibrium response, the kinetics of interfaces, and the nucleation of interfaces independently.
This model was proposed recently in \cite{agrawal2015dynamic}.

For both models, we perform initial-value numerical calculations\footnote
{
    We use large computational domains and take care not to consider results that have been affected by the boundaries. Therefore, we refer to these problems as initial-value problems.
    Strictly speaking, the computational domain is finite and one could refer to them as initial-boundary-value problems.
    However, they aim to mimic a problem on an unbounded domain where the boundaries play no role.
}
to examine the propagation of interfaces.
The results of the initial-value problems, and the analysis in Section \ref{sec:SPFM-discussion}, show that supersonic interfaces cannot be predicted by either of these models.
Therefore, we do not present our attempts to solve traveling wave problems.

%%%%%%%%%%%%%%%%%%%%%
%%%%%%%%%%%%%%%%%%%%%
%%%%%%%%%%%%%%%%%%%%%
%%%%%%%%%%%%%%%%%%%%%
\subsection{Standard Phase-Field Model}
\label{sec:SPFM-formulation}

The standard phase-field model is formulated primarily on the basis of appropriately constructing the energy landscape to obtain the correct equilibrium response.
The evolution is obtained by assuming that $\phi$ follows a steepest descent dynamics, coupled to static or dynamic momentum balance.

The form of the energy of a standard phase-field model coupled with piecewise linear elasticity is:
\begin{equation}
\label{eqn:standard-phase-field-energy}
    F[u,\phi] 
    = 
    \int_\Omega \left(
        \half \alpha |\partial_x\phi|^2 + 
        w(\phi) + \half E(\phi) \left(\partial_x u - \epsilon_0(\phi)\right)^2
    \right) \dm x
\end{equation}
where $w(\phi)$ is a nonconvex energy that favors the formation of interfaces, while the gradient term $\half \alpha |\partial_x\phi|^2$ regularizes them. $\alpha$ is the gradient energy coefficient, and has broadly the same physical interpretation as $\lambda$ in Section \ref{strain gradient model}.

We choose to have phase 1 indicated by $\phi=0$ and phase 2 by $\phi=1$.
Therefore, we choose for $w(\phi)$ the expression:
\begin{equation}
\label{eqn:SPFM-energy-2}
    w(\phi) = \Theta\phi^2(1-\phi)^2 + \phi \Delta\Psi
\end{equation}
The term $\phi^2(1-\phi)^2$ brings in the nonconvexity of the energy landscape, with minima at $\phi=0$ and $\phi=1$.
$\Theta$ is a large constant, chosen to be $10^3$ for this work.
The term $\phi \Delta\Psi$ accounts for the fact that these phases have energy minima at different heights, with the difference quantified by $\Delta\Psi$ that was introduced in \eqref{eqn:stress-energy-response}.
We discuss the relation between the phase-field energy and the classical elastic energy in more detail in Section \ref{sec:app-energies}.

\begin{figure}[htb!]
    %\subfloat[]{\includegraphics[width=0.33\textwidth]{media/W-circ_spfm.jpg}}
    \subfloat[]{\includegraphics[width=0.45\textwidth]{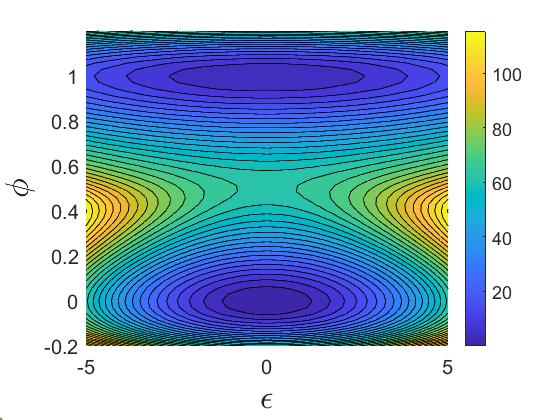}}
    \subfloat[]{\includegraphics[width=0.45\textwidth]{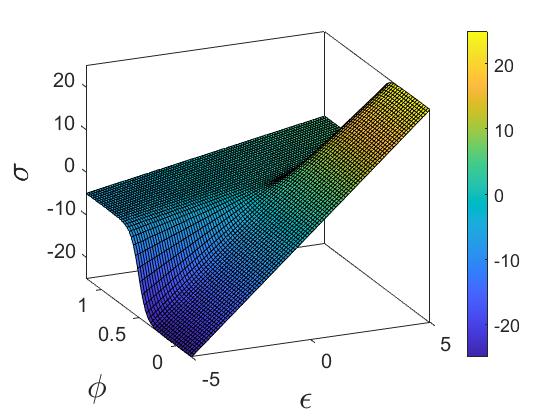}}
    \caption{(a) Plot of the local energy density $w(\phi) + \half E(\phi) \epsilon^2$ as a function of $\epsilon$ and $\phi$. (b) plot of  $\sigma(\epsilon) = E(\phi) \epsilon$ as a function of $\epsilon$ and $\phi$.}
    \label{fig:spfm energy and stress-strain}
\end{figure}

For the elastic response, we set $\epsilon_0(\phi) \equiv 0$, and $E(\phi) = (1 - H_l(\phi-0.5)) E_1 + H_l(\phi-0.5) E_2$, which transitions smoothly from $E(\phi\lesssim 0.5) = E_1$ to $E(\phi\gtrsim0.5) = E_2$.
This mimics the transition between the branches described in \eqref{eqn:stress-energy-response}: we have a linear response with modulus $E_1$ when $\phi\approx 0$ and a linear response with modulus $E_2$ when $\phi\approx 1$.
We notice, however, that the elastic energy $\half E(\phi) (\partial_x u)^2$ is defined for any combination of $\phi$ and $\partial_x u$.
Consequently, it is possible that the system has a strain value corresponding to phase 1 while $\phi\approx 1$; see Figure \ref{fig:spfm energy and stress-strain}.

From \eqref{eqn:standard-phase-field-energy}, we find the stress and driving force:
\begin{align}
\label{eqn:spfm stress and driving force}
    \sigma & = E(\phi) \partial_x u
    \\
    f
    & = - \frac{\delta F}{\delta \phi}  
    = - \left( \parderiv{\ }{\phi} \left(w(\phi) + \half E(\phi)(\partial_x u)^2 \right) - \alpha \partial_{xx} \phi\right)
\end{align}
The evolution equations are given by:
\begin{align}
\label{eqn:spfm evolution eqns}
    & \rho \partial_{tt} u = \partial_x \sigma
    \\
    & \mu \partial_t \phi = f
\end{align}
where we have used the standard steepest descent assumption for $\phi$, and the constant $\mu$ is the mobility.

%%%%%%%%%%%%%%%%%%%%%%%%%%%%%%%%%
%%%%%%%%%%%%%%%%%%%%%%%%%%%%%%%%%
%%%%%%%%%%%%%%%%%%%%%%%%%%%%%%%%%
%%%%%%%%%%%%%%%%%%%%%%%%%%%%%%%%%
\subsubsection{Kinetics of Interfaces from Initial-Value Problems}\label{spfm_kinetics of interface from IVP}

The kinetics of interfaces are studied through numerical solutions of initial-value problems. 
We use standard explicit time-stepping to solve linear momentum balance along with the evolution equation for $\phi$ for the kinetic response.
Our domain is a long finite bar with an interface at the center of the bar. 
Our initial conditions correspond to a displacement / strain field that is not at equilibrium, and the interface has a non-zero driving force across it.
Therefore, the interface will move in the direction of the driving force.  
This is the analog of the classical Riemann problem \cite{Kaushik-peri}.

\begin{figure}[htb!]
    \subfloat[]{\includegraphics[width=0.45\textwidth]{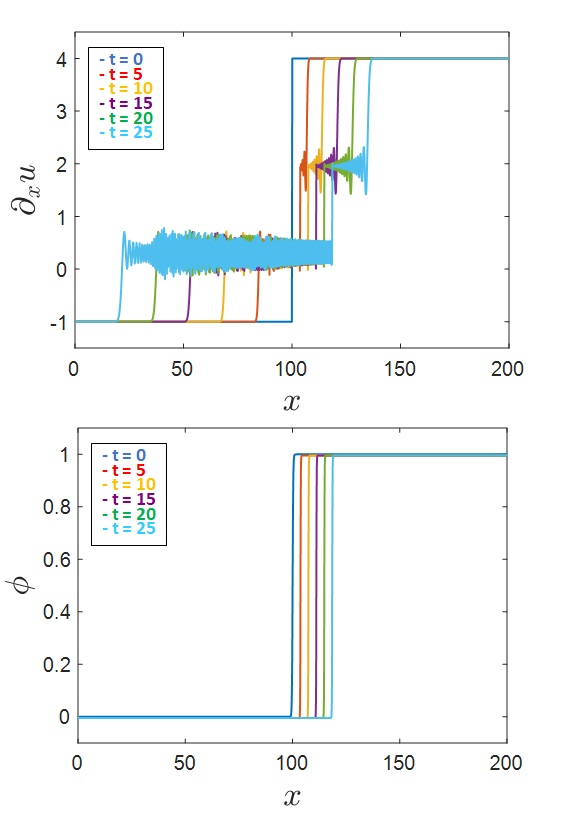}}
    \subfloat[]{\includegraphics[width=0.45\textwidth]{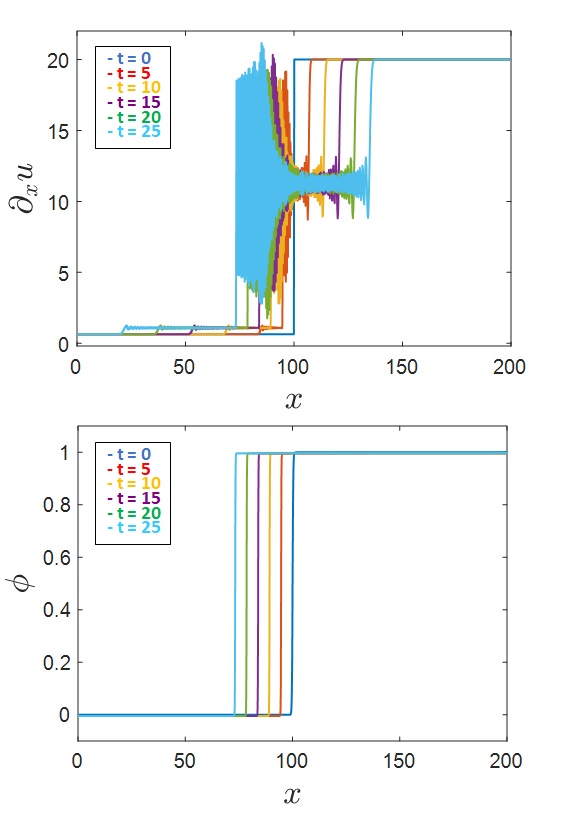}}
    \caption{Initial-value problems for the standard phase-field model with interfaces moving at subsonic velocities. (a) interface moving to the right, into the soft phase; b) interface moving to the left, into the stiff phase. In both cases, we notice acoustic waves propagating in both directions.  The interface can be identified by using that it appears in the plots of both $\phi$ and $\partial_x u$, whereas the acoustic waves have a signature only in $\partial_x u$.}
    \label{fig:Standard phase field model_subsonic}
\end{figure}

Figure \ref{fig:Standard phase field model_subsonic} shows representative results for interfaces that are well below the sonic speed.
These results show the system behaving as we would expect: faster acoustic waves going in both directions from the initial interface, with signatures only in the strain profile; and the slower interface which has a signature in both the strain and $\phi$ profiles.

\begin{figure}[htb!]
    \subfloat[]{\includegraphics[width=0.45\textwidth]{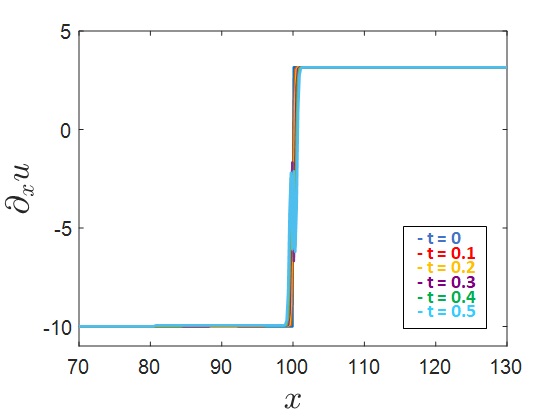}}
    \subfloat[]{\includegraphics[width=0.45\textwidth]{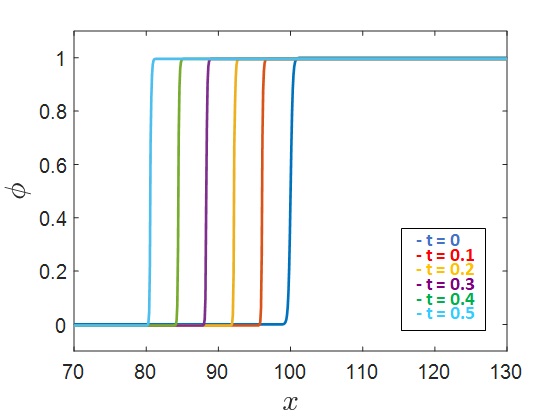}}
    \caption{Initial-value problem for the standard phase-field model close to the sonic velocity. The strain interface has barely moved and is subsonic, while the $\phi$-interface is moving faster than both $c_2$ as well $c_1$!
    While waves above $c_1$ are not permissible by momentum balance, this only constrains the strain and not the evolution of $\phi$.
    Further, the material between the strain interface and the $\phi$ interface is in phase 1 defined through the location of the strain interface but in phase 2 defined through the location of the $\phi$ interface.
    It is therefore in an unphysical part of the energy landscape.}
    \label{fig:Standard phase field model_supersonic}
\end{figure}

An example with a large driving force is shown in Figure \ref{fig:Standard phase field model_supersonic}.
An important feature that we notice is that interfaces in strain and $\phi$ are not at the same location nor do they move with the same velocity.
It is therefore not possible to usefully define an interface velocity.
The strain interface has barely moved and is subsonic, while the $\phi$-interface is moving faster than both $c_2$ as well $c_1$.
Most importantly, the material between the strain interface and the $\phi$ interface is in phase 1 defined through the location of the strain interface but in phase 2 defined through the location of the $\phi$ interface.
It is therefore in an unphysical part of the energy landscape.

%%%%%%%%%%%%%%%%%%%%%%%%%%%%%%%%%
%%%%%%%%%%%%%%%%%%%%%%%%%%%%%%%%%
%%%%%%%%%%%%%%%%%%%%%%%%%%%%%%%%%
%%%%%%%%%%%%%%%%%%%%%%%%%%%%%%%%%
\subsection{Dynamic Phase-field Formulation}
\label{sec:DPFM-formulation}

The dynamic phase-field model aims to transparently separate energetics, kinetics, and nucleation.
In this section, we present the model equations and only those aspects that are directly relevant; we refer to \cite{agrawal2015dynamic,agrawal2015dynamic-2} for the details of the formulation and the characterization.
The kinetics in this model is similar to that proposed earlier by \cite{alber2005solutions} though they have not focused on the energetics or nucleation.

We construct the phase-field energy as follows:
\begin{equation}
\label{eqn:dynamic-phase-field-energy}
    F[u,\phi] = \int_\Omega \left(
        \half \alpha |\partial_x \phi|^2 + (1-H_l(\phi-0.5))\psi_1(\partial_x u) + (H_l(\phi-0.5))\psi_2(\partial_x u)
    \right) \dm x
\end{equation}
where
\begin{equation}
    \psi_A(\partial_x u) 
    =
    W(\epsilon_0^A) + W'(\epsilon_0^A) (\partial_x u - \epsilon_0^A) + \half W''(\epsilon_0^A) (\partial_x u - \epsilon_0^A)^2, \quad A=1,2
\end{equation}
This corresponds to an expansion of the energy about some strain $\epsilon_0^A$ that need not correspond to the stress-free strain.
As a consequence of the piecewise-linearity of the energy, the final expressions turn out to be independent of the choice of $\epsilon_0^1,\epsilon_0^2$; using \eqref{eqn:stress-energy-response} with $\epsilon_0^1 < \epsilon_{1m}$ and $\epsilon_0^2 > \epsilon_{2m}$, we have:
\begin{equation}
    \psi_1(\partial_x u) = \half E_1 (\partial_x u)^2 , \qquad
    \psi_2(\partial_x u) = \half E_2 (\partial_x u)^2 + \Delta\Psi
\end{equation}
We collect expressions and simplify \eqref{eqn:dynamic-phase-field-energy} to get:
\begin{equation}
\label{eqn:dynamic-phase-field-energy-2}
    F[u,\phi] 
    = 
    \int_\Omega \left(
        \half \alpha |\partial_x\phi|^2 + \half E(\phi) \left(\partial_x u \right)^2 + H_l(\phi-0.5) \Delta\Psi
    \right) \dm x
\end{equation}
The key difference is that $w(\phi)$ in \eqref{eqn:standard-phase-field-energy} for the standard energetic phase-field model is replaced by $H_l(\phi-0.5) \Delta\Psi$ above.
We note that $E(\phi) = (1 - H_l(\phi-0.5)) E_1 + H_l(\phi-0.5) E_2$ above, as in the standard phase-field model.

The stress response function and driving force in this phase-field model are:
\begin{align}
    \sigma & = {E(\phi)} \partial_x u
    \\
    f & =
    - \frac{\delta F}{\delta \phi} 
    = - \left( \parderiv{\ }{\phi} \left(H_l(\phi-0.5) \Delta\Psi + \half E(\phi)(\partial_x u)^2 \right) - \alpha \partial_{xx} \phi\right)
    \label{eqn:dpfm stress and driving force}
\end{align}
and the evolution equations are given by:
\begin{align}
    & \rho \partial_{tt} u = \partial_x \sigma
    \\
    \label{eqn:dpfm evolution eqns}
    & \partial_t\phi = \left|\partial_x\phi\right|\hat{v}_n^\phi + G
\end{align}
where $\hat{v}_n^\phi$ is the velocity of the interface and controls the interface kinetics; and $G$ controls the interface nucleation.
In general, both $\hat{v}_n^\phi$ and $G$ can be functions of any quantity such as $f, \sigma, \epsilon$ and their rates; thermodynamics imposes some weak conditions on their dependence on $f$.

We assume the simplest linear kinetics, i.e. $\hat{v}_n^\phi = \kappa f$ with no dependence on other quantities, and $G\equiv 0$.

\begin{figure}[htb!]
    \subfloat[]{\includegraphics[width=0.45\textwidth]{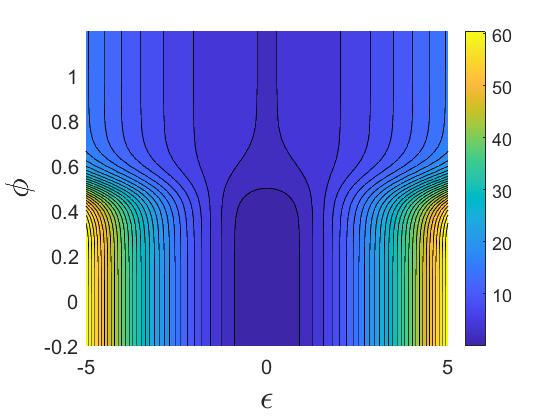}}
    \subfloat[]{\includegraphics[width=0.45\textwidth]{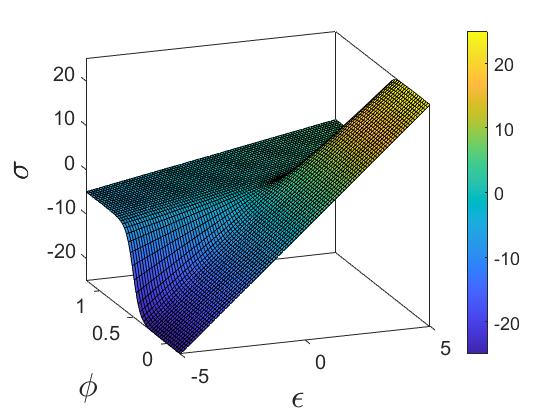}}
    \caption{(a) Plot of the local energy density $H_l(\phi-0.5)\Delta\Psi + \half E(\phi) \epsilon^2$ as a function of $\epsilon$ and $\phi$. (b) plot of  $\sigma(\epsilon) = E(\phi) \epsilon$ as a function of $\epsilon$ and $\phi$.}
    \label{fig:dpfm energy and stress-strain}
\end{figure}

Figure \ref{fig:dpfm energy and stress-strain} plots out the local energy density and the stress.
A key feature of this energy is that it is flat in the $\phi$-direction away from the transition at $\phi \approx 0.5$.
This feature is critical in decoupling nucleation from the structure of the energy landscape.
In particular, this energy landscape simply does not permit nucleation of new phases, regardless of the level of stress / strain.
The only available mechanism for the nucleation of a new phase is through the term denoted $G$ in \eqref{eqn:dpfm evolution eqns}.
This effectively decouples the nucleation of new phases from the equilibrium energetic response, and provides a simple and transparent mechanism to specify the precise conditions for nucleation through the functional dependence of $G$ on $f, \sigma, \partial_x u$, their rates, and so on.

\begin{remark}[Kinetics of Interfaces from Initial-Value Problems]
    As in Section \ref{spfm_kinetics of interface from IVP}, we studied the kinetics of interfaces through numerical solutions of initial-value problems. 
    The results are qualitatively identical to those obtained from the standard phase-field model in Section \ref{spfm_kinetics of interface from IVP}.
\end{remark}

%%%%%%%%%%%%%%%%%%%%%
%%%%%%%%%%%%%%%%%%%%%
%%%%%%%%%%%%%%%%%%%%%
%%%%%%%%%%%%%%%%%%%%%
\section{Unphysical Features of Existing Phase-Field Models}
\label{sec:SPFM-discussion}

We discuss here the 2 main reasons that existing phase-field models are unable to properly handle situations with inertia.
The first reason is related to energetics, namely that the system explores unphysical regions of the energy landscape; further, we highlight that the energy landscape -- even in the energy minimizing setting -- can have unexpected behavior at large strains.
The second reason is that the momentum balance, as typically formulated, leads to a strain singularity for supersonic interfaces.

%%%%%%%%%%%%%%%%%%%%%%%%%%%%%%%%%
%%%%%%%%%%%%%%%%%%%%%%%%%%%%%%%%%
%%%%%%%%%%%%%%%%%%%%%%%%%%%%%%%%%
%%%%%%%%%%%%%%%%%%%%%%%%%%%%%%%%%
\subsection{Energy Landscape at Large Strains}

\begin{figure}[htb!]
    \subfloat[]{\includegraphics[width=0.45\textwidth]{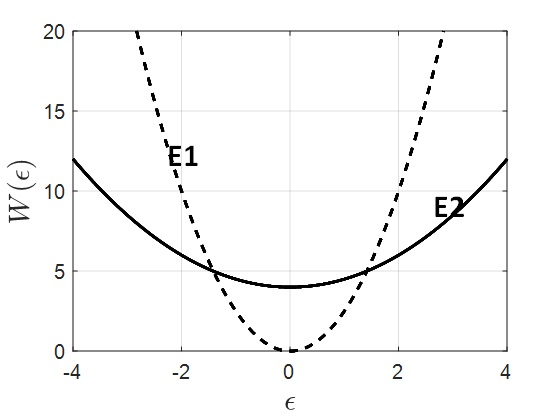}}
    \subfloat[]{\includegraphics[width=0.45\textwidth]{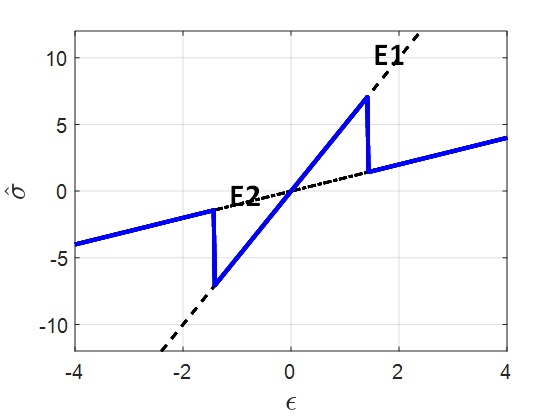}}
    \caption{(a) The energy density for $\phi$ near $0$ and $1$; (b) the stress-strain response for $\phi$ near $0$ and $1$, with the energy minimizing curve highlighted in blue.}
    \label{fig:dpfm_Unphysical behavior at large strains}
\end{figure}

In the 1-d piecewise-linear phase-field energy considered in this paper, the elastic energy has the form $\half E_1 \epsilon^2$ and $\half E_2 \epsilon^2 + \Delta\Psi$ in phases 1 and 2 respectively.
If $E_1>E_2$, we notice that phase 2 always has lower energy when $\epsilon^2$ is large enough, irrespective of the value (or even sign) of $\Delta\Psi$.
This is simply because the quadratic growth eventually wins.
If $\epsilon$ is positive and large, then phase 2 is lower energy; however, even if $\epsilon$ is negative with large magnitude -- considering, e.g., the case of antiplane deformation -- phase 2 again eventually has lower energy.
This unphysical behavior is shown in Figure \ref{fig:dpfm_Unphysical behavior at large strains}.

This conclusion can be readily generalized to higher dimensions and more general energies.
In general, it is a consequence of defining the energy density as a function of $\phi$ and $\epsilon$, and this higher-dimensional energy landscape has several unphysical regions.

In the relatively simple setting considered here, this pathology in the energy landscape effectively introduces an extra, unphysical, transformation that could potentially activate to give unphysical results.
For instance, when we apply large driving forces to attempt to drive the interface at a high velocity, we find that the transformation occurs from the soft phase 2 at high positive strain to the soft phase 2 at high negative strain!
We then observe the unphysical situation of a large strain jump across the interface, but $\phi$ does not evolve at all because the material is in the same phase.

%%%%%%%%%%%%%%%%%%%%%
%%%%%%%%%%%%%%%%%%%%%
%%%%%%%%%%%%%%%%%%%%%
%%%%%%%%%%%%%%%%%%%%%
\subsection{Inability to Model Supersonic Interfaces}\label{Inability to Model Supersonic Interfaces}

As observed in Section \ref{sec:SPFM}, neither of the existing phase-field models could model supersonic interfaces.
We see below that the momentum equation, as used in those models, forbids it.

Consider a steady supersonic interface moving at a given velocity $\dot{s}$.
Using the traveling-wave form $u(x,t) = U(x - \dot{s}t)$ in the momentum equation gives:
\begin{equation}
    \rho\partial_{tt} u = \sigma_x
    \Rightarrow
    \dot{s}^2\rho U'' = \left(E(\phi)U'\right)' 
    \Rightarrow
    \left(\rho \dot{s}^2 - E(\phi)\right) U'= \const
\end{equation}
where we have integrated once to go from the second to the third step above; the constant therefore is independent of $\tilde{x} = x- \dot{s}t$, but can be a function of $\dot{s}$ and material parameters.

Since the interface is supersonic, we have that $c_2 < \dot{s} < c_1$, implying that $E_2 < \rho\dot{s}^2 < E_1$.
Further, we have that $E(\phi)$ varies smoothly between $E_1$ and $E_2$ as $\phi$ transitions from phase 1 to phase 2.
Furthermore, since $\phi(x)$ is smooth, there will be some point $\tilde{x}^*$ in the traveling wave coordinates at which $\rho \dot{s}^2 = E(\phi).$
At $\tilde{x}^*$, we consequently have that $\rho \dot{s}^2 - E(\phi) = 0$, and combining this with $\left(\rho \dot{s}^2 - E(\phi)\right) U'= \const$, we have two possibilities: (1) either $U' \to \infty$ at $\tilde{x}^*$;
or (2) the constant must be $0$, and consequently $U'$ must be zero everywhere except $\tilde{x}^*$.
Neither of these possibilities is acceptable for the interfaces considered here, and therefore we conclude that these phase-field models cannot model supersonic interfaces.

Notice that our argument depends on the continuous variation of $E(\phi)$ with respect to its argument.
If this variation was discontinuous, then we would have to track the moving surface across which $E(\phi)$ was discontinuous, and that would nullify the most important advantage of phase-field modeling, namely that we do not have to track singularities.

%%%%%%%%%%%%%%%%%%%%%
%%%%%%%%%%%%%%%%%%%%%
%%%%%%%%%%%%%%%%%%%%%
%%%%%%%%%%%%%%%%%%%%%
\section{Augmented Dynamic Phase-Field Model for Microstructure Evolution with Inertia}
\label{sec:DPFM}

To address the issues identified in Section \ref{sec:SPFM-discussion}, we propose an augmented phase-field model that has these 2 extra terms:
\begin{enumerate}
    \item a local dynamical term, corresponding to $G$ in \eqref{eqn:dpfm evolution eqns}, that moves the system  -- along the $\phi$ direction only, to avoid disrupting momentum balance -- away from unphysical regions in the energy landscape. This term corresponds to accounting for the missing physics of a driving force that would drive the evolution away from high-/infinite- energy forbidden regions.
    \item a viscous dissipative stress that regularizes the singularities identified in Section \ref{Inability to Model Supersonic Interfaces}, and accounts for the missing physics of dissipative mechanisms that are always active and particularly important at defects and singularities.
\end{enumerate}

%%%%%%%%%%%%%%%%%%%%%
%%%%%%%%%%%%%%%%%%%%%
%%%%%%%%%%%%%%%%%%%%%
%%%%%%%%%%%%%%%%%%%%%

\subsection{Augmented Driving Force to Drive Evolution from Forbidden Regions of the Energy Landscape}\label{nucleation term}

\begin{figure}[htb!]
    \subfloat[]{\includegraphics[width=0.45\textwidth]{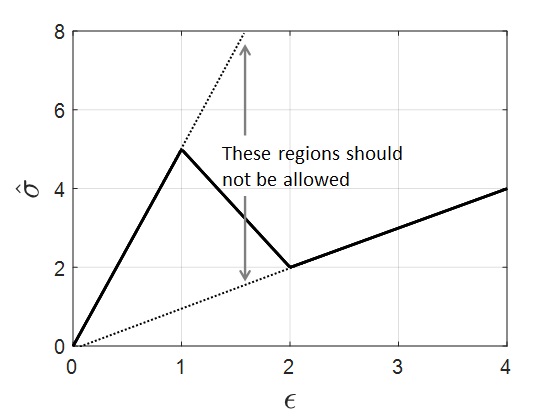}}
    \subfloat[]{\includegraphics[width=0.45\textwidth]{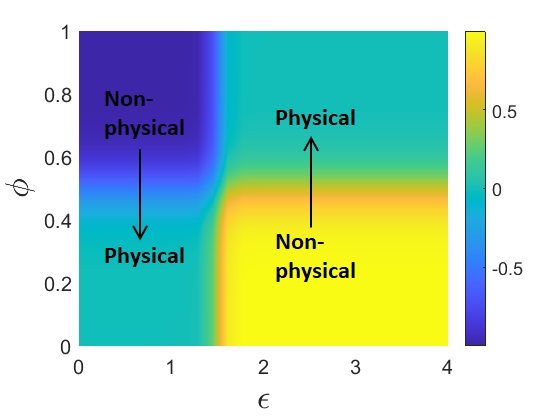}}
    \caption{(a) Disallowed regions in the stress-strain curve. (b) Plot of $G$; the arrows show the evolution of $\phi$ due to $G$.}
    \label{fig:Nucleation term}
\end{figure}

Figure \ref{fig:Nucleation term}(b) shows the energy landscape and the unphysical regions that we aim to avoid.
To avoid introducing new terms in the momentum balance that could lead to spurious unphysical artifacts for shock and acoustic waves, we only modify the evolution of $\phi$ by introducing driving forces in \eqref{eqn:dpfm evolution eqns} through the term denoted by $G$.
In both the standard phase-field model and the dynamic phase-field model, this involves adding a local contribution to the evolution equation for $\phi$, i.e., $\partial_t\phi = \ldots + G$, where $G$ can be a function of any quantity in the system.
$G$ is $0$ when we are in physical regions of the energy landscape; if we enter an unphysical region -- detected by examining the values of $\phi$ and $\epsilon$ -- we set $G$ to be nonzero.
That is, in short, $G$ is a function of $\phi$ and $\epsilon$.
In more general settings, $G$ could be a function of the rates, specific components of the stress or strain, and so on.

Returning to our specific setting, if the system enter the northwest quadrant of Figure \ref{fig:Nucleation term}(b), $G$ is activated and takes on a negative value to push the system downwards into the southwest quadrant.
Similarly, if the system enter the southeast quadrant, $G$ is activated and takes on a positive value to push the system upwards to the northeast quadrant.
The specific expression that we use for $G$ is:
\begin{equation}
\label{eqn:G-expression}
        G(\epsilon,\phi) / G_0
        = 
        H_l\left(\phi - 0.5\right)\left(H_l\left(\epsilon - \bar\epsilon_m\right)-1 \right) 
        + 
        \left(1-H_l\left(\phi - 0.5\right)\right)H_l\left(\epsilon - \bar\epsilon_m\right)
\end{equation}
where $\bar\epsilon_m = \frac{\epsilon_{1m} + \epsilon_{2m}}{2}$, and $G_0$ is a magnitude that we set as large as possible while not losing numerical stability.
The expression above is plotted out in Figure \ref{fig:Nucleation term}(b).

We note that if we interpret $G$ as a nucleation mechanism in the sense of \cite{agrawal2015dynamic}, our choice of $G$ would appear to violate thermodynamic requirements.
However, this apparent violation of thermodynamics is itself an artifact of defining the energy without accounting for the unphysical regions.
A correct definition of the energy landscape would set the energy to $\infty$ in the unphysical regions, but this would be difficult to use for practical computation.
However, in such a landscape where the energy is infinite in the unphysical regions, our choice of $G$ is acceptable to thermodynamics as well as bringing in the missing physics of a driving force to keep the system away from the unphysical regions.

\subsubsection{Quasistatic Characterization}

We characterize the behavior of the phase-field model with the strain/phase constraint -- but without inertia -- to illustrate its effect.
Specifically, we solve numerically the equations:
\begin{align}
    \partial_x \sigma & = 0
    \\
    \mu\partial_t\phi & = \kappa f + G(\partial_x u,\phi)
\end{align}
with load-control and a time-varying applied load.
Because we are in 1-d, this corresponds to simply prescribing $\sigma(t)$; we set this to be a piecewise linear function of time to model loading and unloading.
We start at zero stress and zero strain with $\phi=0$ (phase 1) in the entire specimen; load it until the entire specimen has transformed to phase 2; and then unload back to zero stress and strain.

\begin{figure}[htb!]
    \subfloat[]{\includegraphics[width=0.45\textwidth]{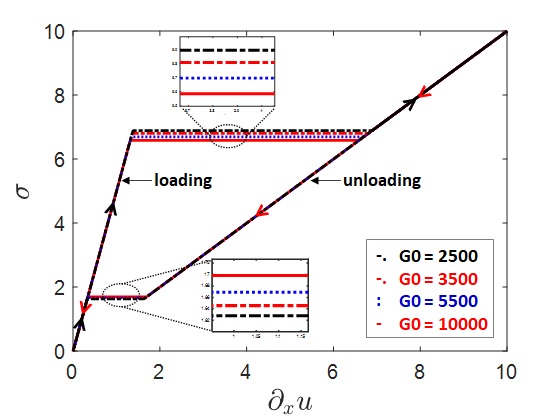}}
    \subfloat[]{\includegraphics[width=0.45\textwidth]{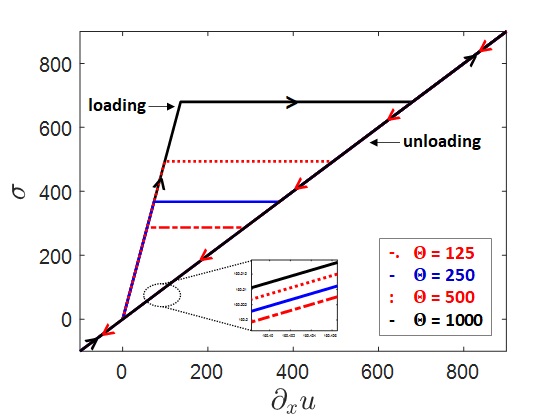}}
    \caption{Stress versus strain under quasistatic loading for the standard phase-field model. (a) With the strain/phase constraint $G$; (b) without $G$. 
    With $G$, we are able to precisely and independently prescribe the critical stress for forward and reverse transformations.
    Without $G$, the system transforms from phase 1 to phase 2 at a critical stress that depends on the energy barrier, that is controlled by $\Theta$. The reverse transformation does not occur at all, even when we go to negative strain; at large negative strains, phase 2 is again stable per the discussion in Section \ref{sec:SPFM-discussion}.}
    \label{fig:spfm_Quasistatic simulation}
\end{figure}

\begin{figure}[htb!]
    \subfloat[]{\includegraphics[width=0.45\textwidth]{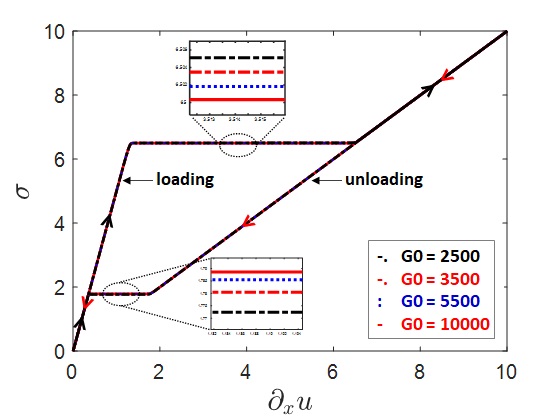}}
    \subfloat[Without $G$, nucleation in both directions is completely suppressed.]{\includegraphics[width=0.45\textwidth]{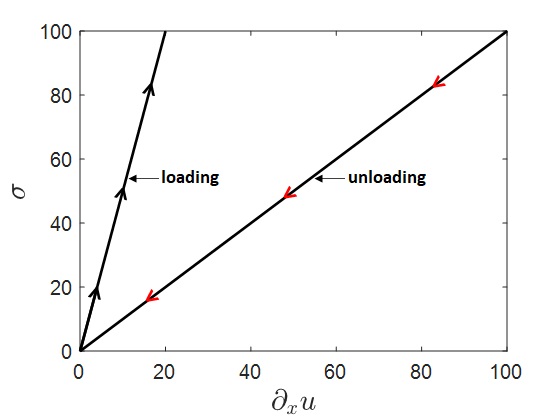}}
    \caption{Stress versus strain under quasistatic loading for the dynamic phase-field model. (a) With the strain/phase constraint $G$; (b) without $G$. 
    With $G$, we are able to precisely and independently prescribe the critical stress for forward and reverse transformations.
    Without $G$, there is simply no transformation in either direction; this is to be expected, as the model was formulated to decouple the nucleation and energetics such that $G$ is the only parameter controlling nucleation \cite{agrawal2015dynamic}.}
    \label{fig:dpfm_Quasistatic simulation}
\end{figure}

The results are shown in Figures \ref{fig:spfm_Quasistatic simulation} and \ref{fig:dpfm_Quasistatic simulation} for the standard and dynamic phase-field models respectively.
For the standard phase-field model, without the strain/phase constraint $G$, we find that the nucleation of the forward transformation is controlled by the height of the energy barrier, while reverse transformation does not occur at all.
On the other hand, with the term $G$ included, nucleation of both forward and reverse transformations occur at precisely the values that we set in the definition of $G$ in \eqref{eqn:G-expression}.
For the dynamic phase-field model, with the term $G$ included, we find similar desirable behavior.
Without $G$, there is no nucleation in either direction, as we desire from the dynamic phase-field model, showing that nucleation is uncoupled and controlled only by the term $G$.

%%%%%%%%%%%%%%%%%%%%%
%%%%%%%%%%%%%%%%%%%%%
%%%%%%%%%%%%%%%%%%%%%
%%%%%%%%%%%%%%%%%%%%%
\subsection{Viscous Dissipative Stresses in the Balance of Momentum}\label{dissipative term}

We recall the simple argument in Section \ref{Inability to Model Supersonic Interfaces} that showed that supersonic interfaces could not be modeled by the existing phase-field models.
In brief, momentum balance in the form $\rho\partial_{tt} u = \partial_x \sigma$, in combination with a traveling wave ansatz, gave $\left(\dot{s}^2\rho - E(\phi)\right)U' = \const$.
Since $\left(\dot{s}^2\rho - E(\phi)\right)$ is $0$ at some location for a supersonic wave, it follows that $U'$ blows up at that point, or is zero essentially everywhere.

If we regularize this equation by using a stress response that includes higher derivatives such as $\partial_{xxt} u$ or $\partial_{xxxx} u$, or indeed any number of other possibilities, we find that our simple argument no longer holds.
We choose to add a linear dissipation of the form $\partial_{xxt} u$, because it is simple to use and rooted in the physics.
The equation of momentum balance will then have the form:
\begin{equation}
    \rho\partial_{tt} u = \partial_x \sigma + \nu\rho\partial_{xxt} u
\end{equation}
We will see in the numerical characterization of the augmented model that this regularization, in combination with the augmented driving force, is sufficient to predict supersonic interfaces.

%%%%%%%%%%%%%%%%%%%%%
%%%%%%%%%%%%%%%%%%%%%
%%%%%%%%%%%%%%%%%%%%%
%%%%%%%%%%%%%%%%%%%%%
\subsection{Characterization of Augmented Phase-Field Models}

%%%%%%%%%%%%%%%%%%%%%
%%%%%%%%%%%%%%%%%%%%%
%%%%%%%%%%%%%%%%%%%%%
%%%%%%%%%%%%%%%%%%%%%

\subsubsection{Augmented Standard Phase-Field Model}
\label{sec:augmented-SPFM-results}

Using the model from Section \ref{sec:SPFM-formulation}, we have the following expressions for the free energy, stress, and driving force:
\begin{equation}
\begin{split}
    &F[u,\phi] = \int_\Omega \left(
        \half \alpha |\partial_x\phi|^2 + w(\phi) + \half E(\phi) \left(\partial_x u \right)^2
    \right) \dm x
    \\
    & \Rightarrow
    \sigma = E(\phi) \partial_x u
    ,\quad
    f 
    = - \left( \parderiv{\ }{\phi} \left(w(\phi) + \half E(\phi)\left(\partial_x u \right)^2 \right) - \alpha \partial_{xx} \phi\right)
\end{split}
\end{equation}

With the viscous stress and the augmented driving force, the evolution equations are given by:
\begin{align}
\label{eqn:augmented-SPFM-evolution}
    \rho \partial_{tt} u & = \partial_x \sigma + \nu\rho\partial_{xxt} u
    \\
    \mu \partial_t \phi & = f + G(\partial_x u,\phi)
\end{align}

\begin{figure}[htb!]
    \subfloat[]{\includegraphics[width=0.45\textwidth]{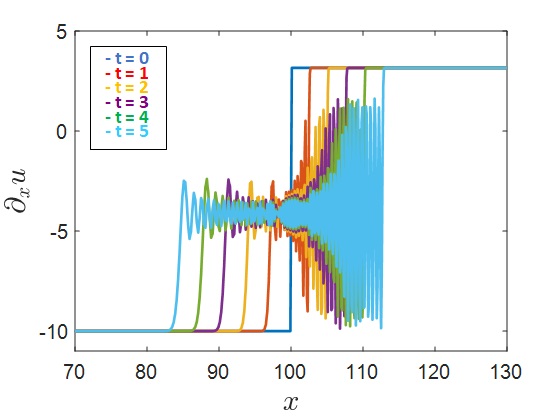}}
    \subfloat[]{\includegraphics[width=0.45\textwidth]{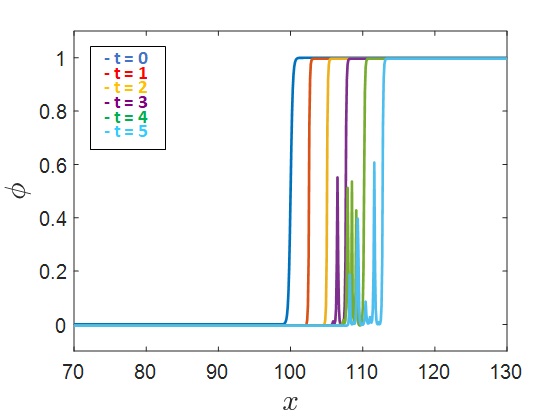}}
    \caption{Initial-value problems for the standard phase-field model with supersonic interfaces, using only the augmented driving force $G$. (a) $\partial_x u$, (b) $\phi$. We see that $\phi$ has undesirable oscillations, but the interfaces in strain and $\phi$ move together, showing that the system is not in unphysical regions of the energy landscape.}
    \label{fig:Standard phase field model _supersonic_G}
\end{figure}

\begin{figure}[htb!]
    \subfloat[]{\includegraphics[width=0.45\textwidth]{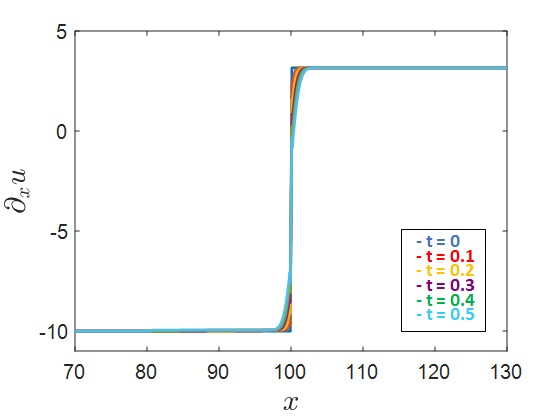}}
    \subfloat[]{\includegraphics[width=0.45\textwidth]{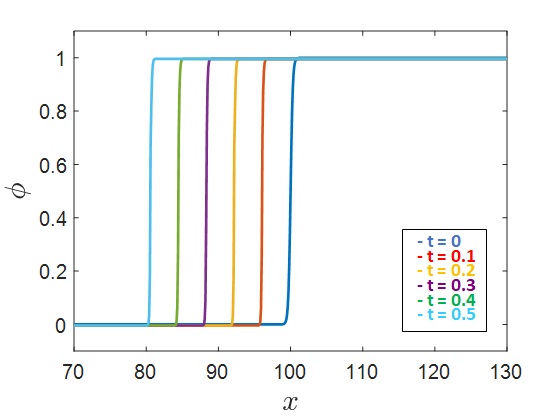}}
    \caption{Initial-value problems for the standard phase-field model with supersonic interfaces, using only viscous stresses. (a) $\partial_x u$, (b) $\phi$. We see, as before, that the strain interface is subsonic while the $\phi$ interface is well above supersonic. Therefore, the system explores unphysical regions of the energy landscape.}
    \label{fig:Standard phase field model _supersonic_Damping}
\end{figure}

\begin{figure}[htb!]
    \subfloat[]{\includegraphics[width=0.45\textwidth]{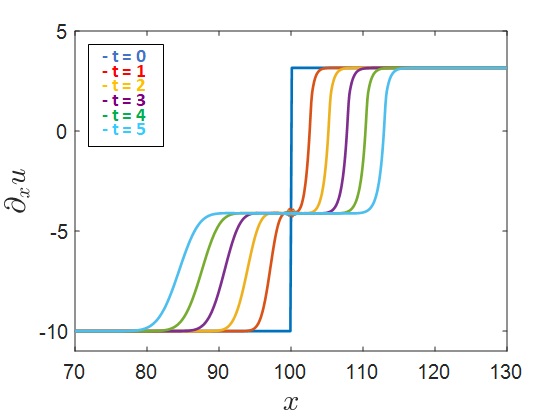}}
    \subfloat[]{\includegraphics[width=0.45\textwidth]{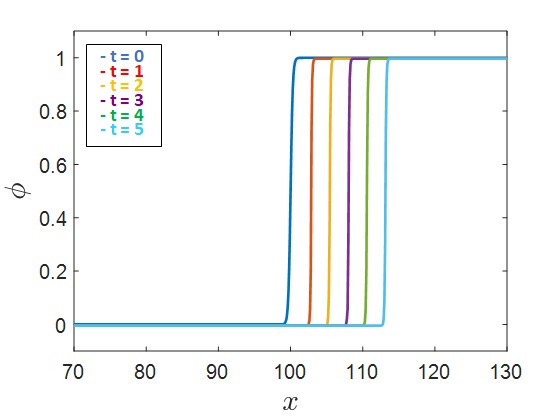}}
    \caption{Initial-value problems for the standard phase-field model with supersonic interfaces, using both the augmented driving force and viscous stresses. (a) $\partial_x u$, (b) $\phi$. The evolution is precisely as we desire, in that the strain and $\phi$ interfaces move together supersonically with no undesirable oscillations. Notice that the strain interfaces are smeared out due to the additional dissipative regularization. Note that the ``blip'' in $\partial_x u$ at $t=1$ is a transient that has not yet stabilized into a steadily-moving interface.}
    \label{fig:Standard phase field model _supersonic_G+Damping}
\end{figure}

Numerical solutions of initial-value problems with supersonic interfaces are shown in Figures \ref{fig:Standard phase field model _supersonic_G}, \ref{fig:Standard phase field model _supersonic_Damping}, and \ref{fig:Standard phase field model _supersonic_G+Damping}, for the case with only the augmented driving force, the case with only viscous stress, and the case with both mechanisms respectively.
We present the results with only the augmented driving force and only the viscous stress to show that both of these mechanisms are essential.

\begin{figure}[htb!]
    {\includegraphics[width=0.6\textwidth]{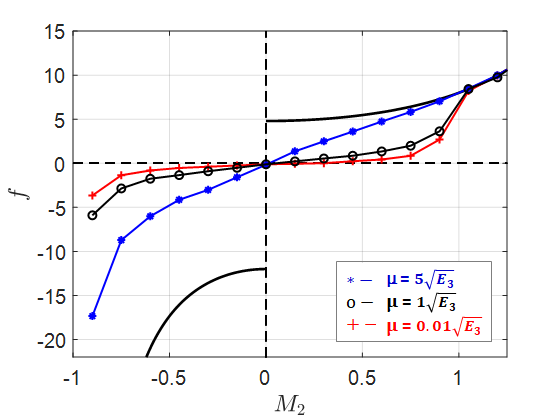}}
    \caption{The induced kinetic relation in the augmented standard phase-field model with dissipative regularization and strain/phase constraints, for different choices of model parameters, computed using traveling wave calculations.}
    \label{fig:spfm+G+damp-kinetics}
\end{figure}

\begin{figure}[htb!]
    \subfloat[]{\includegraphics[width=0.33\textwidth]{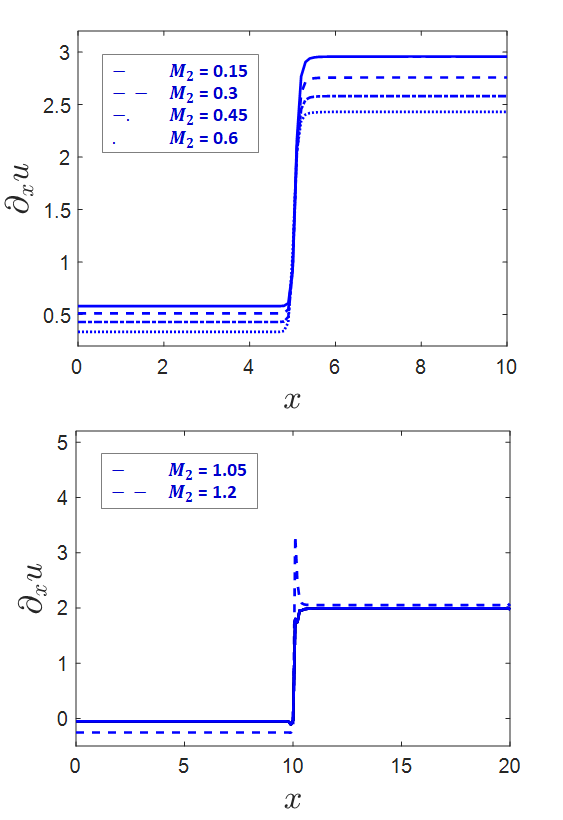}}
    \subfloat[]{\includegraphics[width=0.33\textwidth]{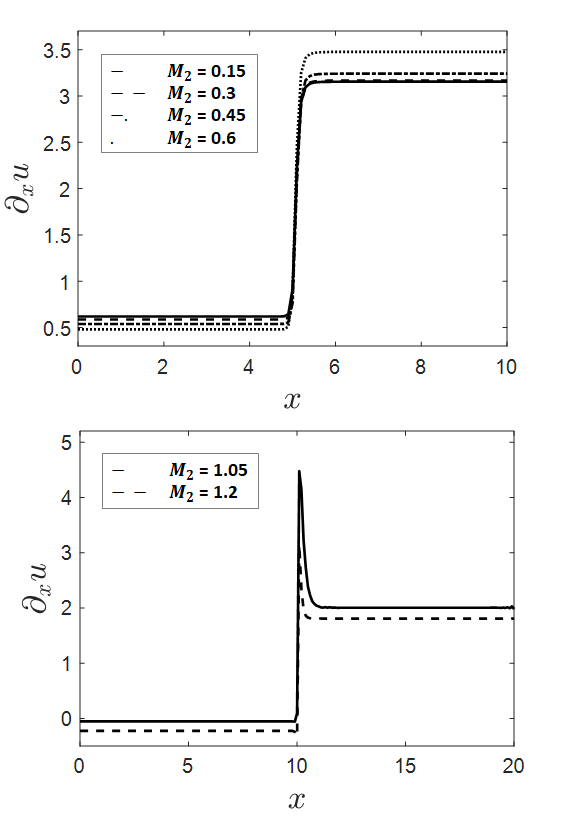}}
    \subfloat[]{\includegraphics[width=0.33\textwidth]{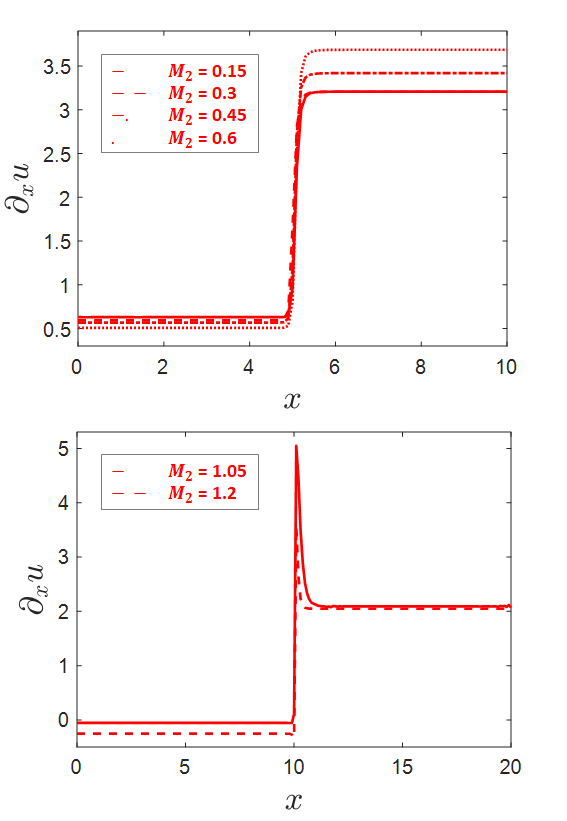}}
    \caption{Representative strain profiles from the traveling wave computations with the standard phase-field model with dissipative regularization and strain/phase constraints. Parameters : $\nu=0.02, G_0=1, \alpha=10$. Representative traveling wave strain profiles for different velocities and different values of $\mu$. (a) $\mu=5\sqrt{\frac{E_3}{\rho}}$, (b) $\mu=1\sqrt{\frac{E_3}{\rho}}$, (c) $\mu=0.01\sqrt{\frac{E_3}{\rho}}$.}
    \label{fig:spfm+G+damp-profiles}
\end{figure}

We next perform traveling wave calculations to find the kinetic relations for various choices of model parameters.
The overall numerical approach is similar to that described in Section \ref{strain gradient model_traveling waves}.
The primary difference is that we have 2 simultaneous equations to solve in \eqref{eqn:augmented-SPFM-evolution}.
Our solution strategy is to compute the residual for each equation, square both residuals, and then add the squared residuals and integrate over the domain to obtain a single functional that we can minimize.
Figure \ref{fig:spfm+G+damp-kinetics} shows that the kinetics is sensitive to the model parameters for subsonic interfaces, but is not for supersonic interfaces; Figure \ref{fig:spfm+G+damp-profiles} shows some representative strain profiles.

%%%%%%%%%%%%%%%%%%%%%
%%%%%%%%%%%%%%%%%%%%%
%%%%%%%%%%%%%%%%%%%%%
%%%%%%%%%%%%%%%%%%%%%

\subsubsection{Augmented Dynamic phase-field model}
\label{sec:augmented-DPFM-results}

Using the model from Section \ref{sec:DPFM-formulation}, we have the following expressions for the free energy, stress, and driving force:
\begin{equation}
\begin{split}
    &F[u,\phi] = 
    \int_\Omega \left(
        \half \alpha |\partial_x\phi|^2 + \half E(\phi) \left(\partial_x u \right)^2 + H_l(\phi-0.5) \Delta\Psi
    \right) \dm x
    \\
    & \Rightarrow
    \sigma = {E(\phi)} \partial_x u,
    \quad    
    f
    = 
     - \left( \parderiv{\ }{\phi} \left(H_l(\phi-0.5) \Delta\Psi + \half E(\phi)(\partial_x u)^2 \right) - \alpha \partial_{xx} \phi\right)
\end{split}
\end{equation}

With both the viscous stress and the augmented driving force included, the evolution equations are given by:
\begin{align}
\label{eqn:augmented-DPFM-evolution}
    & \rho \partial_{tt} u = \partial_x \sigma + \nu\rho\partial_{xxt} u
    \\
    & \partial_t\phi = \left|\partial_x\phi\right|\hat{v}_n^\phi+ G(\partial_x u,\phi)
\end{align}
and we use linear kinetics $\hat{v}_n^\phi = \kappa f$ as in Section \ref{sec:DPFM-formulation}.
We perform numerical computations of both initial-value problems and traveling wave problems.
The results of the initial-value problems are qualitatively identical to those obtained with the augmented standard phase-field model reported in Section \ref{sec:augmented-SPFM-results}, and we do not present the details.
In summary, neither viscous stresses nor the augmented driving force by themselves lead to good results -- and the bad results are qualitatively similar to those reported in Section \ref{sec:augmented-SPFM-results} -- while the use of both mechanisms together provide the desired results.

\begin{figure}[htb!]
    {\includegraphics[width=0.6\textwidth]{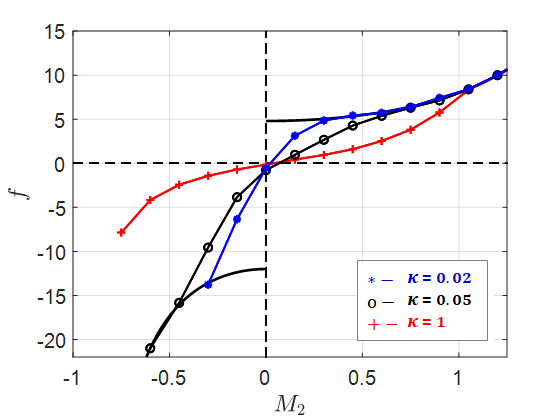}}
    \caption{The induced kinetic relation in the augmented dynamic phase-field model with dissipative regularization and strain/phase constraints, for different choices of model parameters, computed using traveling wave calculations.}
    \label{fig:dpfm+G+damp-kinetics}
\end{figure}

\begin{figure}[htb!]
    \subfloat[]{\includegraphics[width=0.33\textwidth]{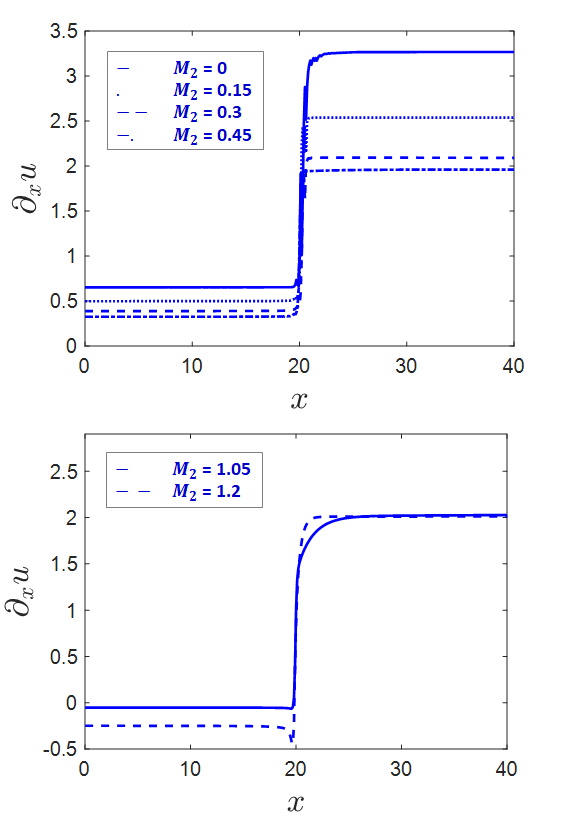}}
    \subfloat[]{\includegraphics[width=0.33\textwidth]{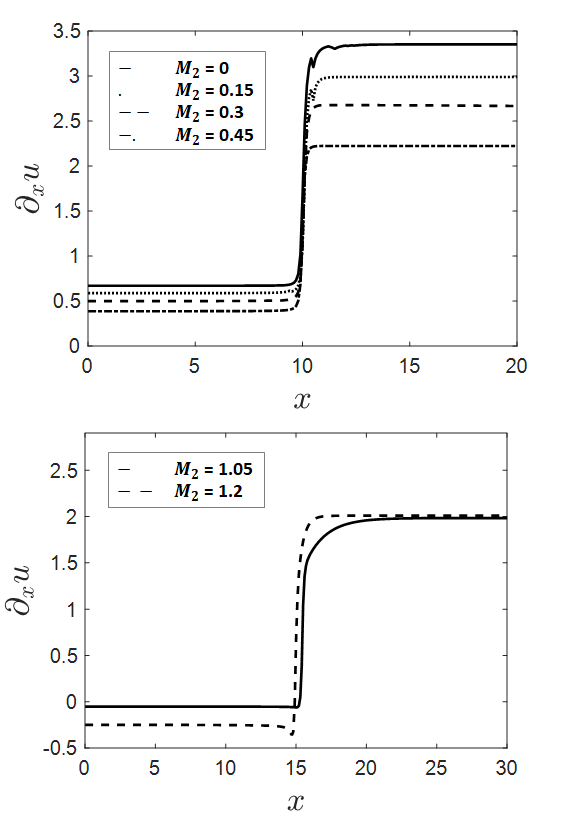}}
    \subfloat[]{\includegraphics[width=0.33\textwidth]{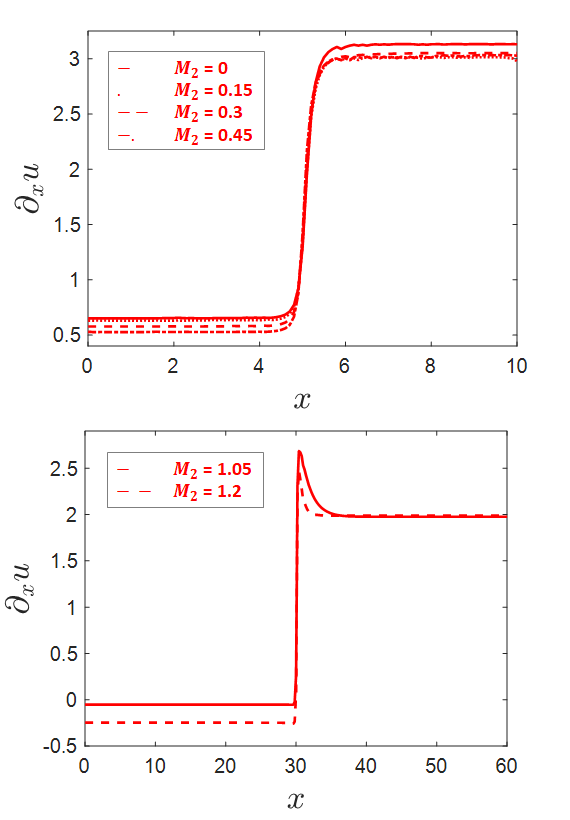}}
    \caption{Representative strain profiles from the traveling wave computations with the dynamic phase-field model with dissipative regularization and strain/phase constraints. Parameters : $\nu=0.2, G_0=1, \alpha=20$. Representative traveling wave strain profiles for different velocities and different values of $\kappa$. (a) $\kappa=0.02$, (b) $\kappa=0.05$, (c) $\kappa=1$.}
    \label{fig:dpfm+G+damp-profiles}
\end{figure}

The traveling wave calculations of interface kinetics also show the desired behavior.
Namely, the kinetic relation for subsonic interfaces is sensitive to model parameters, while the kinetics of supersonic interfaces is not (Figure \ref{fig:dpfm+G+damp-kinetics}); representative strain profiles are shown in Figure \ref{fig:dpfm+G+damp-profiles}.

%%%%%%%%%%%%%%%%%%%%%
%%%%%%%%%%%%%%%%%%%%%
%%%%%%%%%%%%%%%%%%%%%
%%%%%%%%%%%%%%%%%%%%%
\section{Discussion}
\label{sec:discussion}

Dynamic interfaces in a non-monotone stress-strain material are predicted by classical elastodynamics to have two regimes: subsonic where the evolution is nonunique, and supersonic where the evolution is unique\footnote{
    An interesting analogy to this appears in the work of Eshelby on liquid crystals, where he found that the configurational force on a disclination is identical to the real force \cite{eshelby1980force,cermelli2002evolution,singh2020pseudomomentum}.
    The coexistence of regimes with uniqueness and nonuniqueness within a given problem have been explored also in the context of soft materials, e.g. \cite{cohen2019dynamic,knowles2002impact}.
}.
We show that peridynamics preserves this key feature of interfaces (Figure \ref{fig:peridynamics-kinetics}), but existing phase-field model do not; in fact, supersonic interfaces are not admitted by existing phase-field models.
We propose an augmentation of phase-field models, using a viscous stress in the momentum balance and a local dynamical term that keeps the system out of unphysical regions of the energy landscape.
We demonstrate that the augmented phase-field models recover the feature that subsonic interfaces are sensitive to model parameters while supersonic interfaces are not (Figures \ref{fig:spfm+G+damp-kinetics} and \ref{fig:dpfm+G+damp-kinetics}).

A contribution of this paper is in providing a critical qualitative test that distinguishes between the predictions of peridynamics and phase-field models.
Prior work has shown that both peridynamics and phase-field models trivially recover the homogeneous deformation limit, and rigorous work for both peridynamics and phase-field models show that they recover the energy-minimizing Griffith theory of brittle fracture without inertia \cite{lipton2014dynamic,ambrosio1990approximation}.
Recent numerical works compare against experiment \cite{mehrmashhadi2020validating,dally2020cohesive}; while valuable and complementary, these leave open the question if the failure to reproduce an experiment is a calibration issue or a difference in the fundamental structure of the models.
This paper shows that there are fundamental differences in the dynamic setting that cannot be bridged by calibration of model parameters.

A key shortcoming of existing phase-field models that is discussed in this paper is that the energy landscape is expanded by introducing the phase-field $\phi$, and this expanded energy landscape has numerous unphysical regions.
For equilibrium problems, the system is governed by energy minimization, and hence it does not explore these unphysical regions.
However, in dynamic problems with inertia, the potential energy is balanced against the kinetic energy.
Therefore, energy-minimization formulations of the energy are found to be inadequate to avoid the unphysical regions, and we require additional physics to prevent the system from exploring the unphysical regions.
In short, if new equations and variables are introduced, we require additional physics to ensure that they behave appropriately in all regimes of application.
Our proposed augmentation of phase-field models includes this additional physics in the form of an additional driving force.

We have examined in detail two types of phase-field models, but we expect similar results from other phase-field models, e.g. \cite{kamensky2018hyperbolic}, because the key features that lead to these findings are similar in those other models.
In turn, we expect that the augmentation proposed in this paper will prove useful in augmenting also those other models.
It will similarly be interesting to examine if other regularized models for interfaces, e.g. \cite{clayton2017finsler,clayton2019nonlinear,dayal2017leading}, show similar results.
Related to this, an important next step is to test the augmented model in realistic higher-dimensional problems, such as dynamic fracture.
For instance, a key test is to perform numerical calculations of dynamic fracture with an augmented phase-field model and compare it quantitatively to the predictions from an existing phase-field model.
In addition to the augmentation directly affecting the crack growth dynamics, it will also affect the kinetics of nonlinear waves that govern the transport of elastic energy around growing defects \cite{marder2006supersonic,gao1997elastic,cohen2015dynamic,kumar2017some,faye2017spherical}.

%%%%%%%%%%%%%%%%%%%%%
%%%%%%%%%%%%%%%%%%%%%
%%%%%%%%%%%%%%%%%%%%%
%%%%%%%%%%%%%%%%%%%%%

\section*{Software Availability}

A version of the code developed for this work is available at  \\ \url{https://github.com/janelchua/Phase-field-and-Peridynamics.git}

\begin{acknowledgments}
    We thank G. Ravichandran, Anthony Rollett, and Phoebus Rosakis for useful discussions and encouragement;
    Army Research Office (MURI W911NF-19-1-0245, W911NF-17-1-0084), Office of Naval Research (N00014-18-1-2528), and National Science Foundation (CMMI MOMS 1635407, DMREF 2118945, DMS 2108784) for support; and
    National Science Foundation for XSEDE resources provided by Pittsburgh Supercomputing Center.
\end{acknowledgments}

%%%%%%%%%%%%%%%%%%%%%
%%%%%%%%%%%%%%%%%%%%%
%%%%%%%%%%%%%%%%%%%%%
%%%%%%%%%%%%%%%%%%%%%

\appendix

\section{Correspondence between Strain Energy Densities of Classical Elasticity and Phase-Field Models}
\label{sec:app-energies}

The Maxwell stress is an important physical quantity that characterizes phase transformations modeled by nonconvex energies.
It is the value of the mechanical stress at which the driving force on the interface \eqref{eqn:classical-driving-force} is $0$ in the quasistatic setting.
Because the driving force for interface motion is $0$, the (infinitesimal) motion of the interface in either direction does not change the energy of the system.
Therefore, equivalently, it is the value of the mechanical stress at which there is no energetic preference for either phase \cite{abeyaratne2006evolution}.
It is essential that the energetics of all the phase-field models proposed here give rise to the same Maxwell stress as predicted by classical elasticity.
We show here that our construction of these models satisfies this requirement.

First, we compute the Maxwell stress for the classical elasticity model in \eqref{eqn:stress-energy-response}.
Let the driving force $f$ vanish at the Maxwell stress $\sigma_M$.
From \eqref{eqn:classical-driving-force}, we get the condition:
\begin{equation}
\label{eqn:Maxwell}
    \llbracket W(\epsilon) - \sigma_M \epsilon \rrbracket = 0
\end{equation}
Since we have assumed quasistatics, we have that $\sigma_M = E_1 \epsilon^- = E_2 \epsilon^+$, where $\epsilon^\pm$ are the uniform strains on either side of the interface.
Further, using \eqref{eqn:stress-energy-response}, we have that $W(\epsilon^+) = \half E_2 \epsilon^{+2} + \Delta\Psi$ and $W(\epsilon^-) = \half E_1 \epsilon^{-2}$.
Substituting these in \eqref{eqn:Maxwell}, we find the expression for the Maxwell stress:
\begin{equation}
    \sigma_M = \sqrt{\frac{2 \Delta\Psi}{\frac{1}{E_2}-\frac{1}{E_1}}}
\end{equation}

We now consider the energetics of the standard phase-field model discussed in Section \ref{sec:SPFM}.
The local part of the energy density is $w(\phi) + \half E(\phi) \epsilon^2$.
Decompose $w(\phi) = \tilde{w}(\phi) + \phi \Delta\Psi$.
The nonconvex energy $\tilde{w}(\phi)$ is given a specific form in that section, but for our discussion we only need that it contains two minima ($\phi = 0, 1$) and that these minima are at the same height.
We notice that the energy difference between the the two phases at zero stress is $\Delta\Psi$.

Next we consider the Maxwell stress predicted by the standard phase-field model from  \eqref{eqn:standard-phase-field-energy}, \eqref{eqn:SPFM-energy-2}.
Using the interpretation of the Maxwell stress as the stress at which the {\em potential} energy difference between the phases is zero gives:
\begin{equation}
    \half E_2 \epsilon^{-2} + \Delta\Psi -\sigma_M \epsilon^- = \half E_1 \epsilon^{+2} - \sigma_M \epsilon^+
    \Rightarrow
    \sigma_M = \sqrt{\frac{2 \Delta\Psi}{\frac{1}{E_2}-\frac{1}{E_1}}}
\end{equation}
where we have used $\sigma_M = E_1 \epsilon^+ = E_2 \epsilon^-$, and $\epsilon^{\pm}$ are the uniform far-field strains.

We highlight an important approximation in our calculation above.
Under stress, the energy minima are not precisely at $\phi=0,1$, and hence the difference in energy between the phases is not precisely $\Delta\Psi$.
This calculation is exact only in the limit that $\Theta$ is large.

We next examine the energy of the dynamic phase-field model from \eqref{eqn:dynamic-phase-field-energy-2}.
The local part of the energy is $H_l(\phi-0.5)\Delta\Psi + \half E(\phi) \epsilon^2$.
The difference between the energies at zero strain is $\Delta\Psi$.
Therefore, the value of the Maxwell stress in this model is also identical to the values in the standard energetic phase-field and classical elasticity.
Further, we notice that this result is exact because the difference between the energies is independent of stress, i.e., $H_l(\phi-0.5)$ goes to $0$ and $1$ regardless of the stress.

%%%%%%%%%%%%%%%%%%%%%
%%%%%%%%%%%%%%%%%%%%%
%%%%%%%%%%%%%%%%%%%%%
%%%%%%%%%%%%%%%%%%%%%

% References
\bibliographystyle{alpha}
\bibliography{kinetics-refs}

\newcommand{\etalchar}[1]{$^{#1}$}
\begin{thebibliography}{LMTS{\etalchar{+}}18}

\bibitem[AA12]{abdollahi_arias}
A.~Abdollahi and I.~Arias.
\newblock Phase-field modeling of crack propagation in piezoelectric and
  ferroelectric materials with different electromechanical crack conditions.
\newblock {\em Journal of the Mechanics and Physics of Solids}, 60:2100--2126,
  2012.

\bibitem[AD15a]{agrawal2015dynamic}
Vaibhav Agrawal and Kaushik Dayal.
\newblock A dynamic phase-field model for structural transformations and
  twinning: Regularized interfaces with transparent prescription of complex
  kinetics and nucleation. part i: Formulation and one-dimensional
  characterization.
\newblock {\em Journal of the Mechanics and Physics of Solids}, 85:270--290,
  2015.

\bibitem[AD15b]{agrawal2015dynamic-2}
Vaibhav Agrawal and Kaushik Dayal.
\newblock A dynamic phase-field model for structural transformations and
  twinning: Regularized interfaces with transparent prescription of complex
  kinetics and nucleation. part ii: Two-dimensional characterization and
  boundary kinetics.
\newblock {\em Journal of the Mechanics and Physics of Solids}, 85:291--307,
  2015.

\bibitem[AD17]{agrawal2017dependence}
Vaibhav Agrawal and Kaushik Dayal.
\newblock Dependence of equilibrium griffith surface energy on crack speed in
  phase-field models for fracture coupled to elastodynamics.
\newblock {\em International Journal of Fracture}, 207(2):243--249, 2017.

\bibitem[AGDL15]{ambati2015review}
Marreddy Ambati, Tymofiy Gerasimov, and Laura De~Lorenzis.
\newblock A review on phase-field models of brittle fracture and a new fast
  hybrid formulation.
\newblock {\em Computational Mechanics}, 55(2):383--405, 2015.

\bibitem[AHKB20]{albrecht2020phase}
C~Albrecht, A~Hunter, A~Kumar, and IJ~Beyerlein.
\newblock A phase field model for dislocations in hexagonal close packed
  crystals.
\newblock {\em Journal of the Mechanics and Physics of Solids}, 137:103823,
  2020.

\bibitem[AK90]{abeyaratne1990driving}
Rohan Abeyaratne and James~K Knowles.
\newblock On the driving traction acting on a surface of strain discontinuity
  in a continuum.
\newblock {\em Journal of the Mechanics and Physics of Solids}, 38(3):345--360,
  1990.

\bibitem[AK91a]{abeyaratne1991implications}
Rohan Abeyaratne and James~K Knowles.
\newblock Implications of viscosity and strain-gradient effects for the
  kinetics of propagating phase boundaries in solids.
\newblock {\em SIAM Journal on Applied Mathematics}, 51(5):1205--1221, 1991.

\bibitem[AK91b]{abeyaratne1991kinetic}
Rohan Abeyaratne and James~K Knowles.
\newblock Kinetic relations and the propagation of phase boundaries in solids.
\newblock {\em Archive for rational mechanics and analysis}, 114(2):119--154,
  1991.

\bibitem[AK06]{abeyaratne2006evolution}
Rohan Abeyaratne and James~K Knowles.
\newblock {\em Evolution of phase transitions: a continuum theory}.
\newblock Cambridge University Press, 2006.

\bibitem[AT90]{ambrosio1990approximation}
Luigi Ambrosio and Vincenzo~Maria Tortorelli.
\newblock Approximation of functional depending on jumps by elliptic functional
  via $gamma$-convergence.
\newblock {\em Communications on Pure and Applied Mathematics},
  43(8):999--1036, 1990.

\bibitem[AZ05]{alber2005solutions}
Hans-Dieter Alber and Peicheng Zhu.
\newblock Solutions to a model with nonuniformly parabolic terms for phase
  evolution driven by configurational forces.
\newblock {\em SIAM Journal on Applied Mathematics}, 66(2):680--699, 2005.

\bibitem[BD18]{breitzman2018bond}
Timothy Breitzman and Kaushik Dayal.
\newblock Bond-level deformation gradients and energy averaging in
  peridynamics.
\newblock {\em Journal of the Mechanics and Physics of Solids}, 110:192--204,
  2018.

\bibitem[BH16]{beyerlein2016understanding}
IJ~Beyerlein and A~Hunter.
\newblock Understanding dislocation mechanics at the mesoscale using phase
  field dislocation dynamics.
\newblock {\em Philosophical Transactions of the Royal Society A: Mathematical,
  Physical and Engineering Sciences}, 374(2066):20150166, 2016.

\bibitem[BRLM17]{bleyer2017dynamic}
J{\'e}r{\'e}my Bleyer, Cl{\'e}ment Roux-Langlois, and Jean-Fran{\c{c}}ois
  Molinari.
\newblock Dynamic crack propagation with a variational phase-field model:
  limiting speed, crack branching and velocity-toughening mechanisms.
\newblock {\em International Journal of Fracture}, 204(1):79--100, 2017.

\bibitem[BVS{\etalchar{+}}12]{borden2012phase}
Michael~J Borden, Clemens~V Verhoosel, Michael~A Scott, Thomas~JR Hughes, and
  Chad~M Landis.
\newblock A phase-field description of dynamic brittle fracture.
\newblock {\em Computer Methods in Applied Mechanics and Engineering},
  217:77--95, 2012.

\bibitem[CF02]{cermelli2002evolution}
Paolo Cermelli and Eliot Fried.
\newblock The evolution equation for a disclination in a nematic liquid
  crystal.
\newblock {\em Proceedings of the Royal Society of London. Series A:
  Mathematical, Physical and Engineering Sciences}, 458(2017):1--20, 2002.

\bibitem[Che02]{Chen_Phasefield}
Long-Qing Chen.
\newblock Phase-field models for microstructure evolution.
\newblock {\em Annual Review of Materials Research}, 32:113--140, 2002.

\bibitem[CK14]{clayton2014geometrically}
John~D Clayton and J~Knap.
\newblock A geometrically nonlinear phase field theory of brittle fracture.
\newblock {\em International Journal of Fracture}, 189(2):139--148, 2014.

\bibitem[Cla17]{clayton2017finsler}
JD~Clayton.
\newblock Finsler-geometric continuum dynamics and shock compression.
\newblock {\em International Journal of Fracture}, 208(1-2):53--78, 2017.

\bibitem[Cla19]{clayton2019nonlinear}
John~D Clayton.
\newblock {\em Nonlinear elastic and inelastic models for shock compression of
  crystalline solids}.
\newblock Springer, 2019.

\bibitem[CM15]{cohen2015dynamic}
Tal Cohen and Alain Molinari.
\newblock Dynamic cavitation and relaxation in incompressible nonlinear
  viscoelastic solids.
\newblock {\em International Journal of Solids and Structures}, 69:544--552,
  2015.

\bibitem[CMB20]{collet2020variational}
Sylvain Collet, Jean-Fran{\c{c}}ois Molinari, and Stella Brach.
\newblock Variational phase-field continuum model uncovers adhesive wear
  mechanisms in asperity junctions.
\newblock {\em Journal of the Mechanics and Physics of Solids}, 145:104130,
  2020.

\bibitem[Coh19]{cohen2019dynamic}
Tal Cohen.
\newblock Dynamic enlargement of a hole in a sheet: Crater formation and
  propagation of cylindrical shock waves.
\newblock {\em Journal of the Mechanics and Physics of Solids}, 133:103743,
  2019.

\bibitem[Daf05]{dafermos2005hyperbolic}
Constantine~M Dafermos.
\newblock {\em Hyperbolic conservation laws in continuum physics}.
\newblock Springer, 2005.

\bibitem[Day17]{dayal2017leading}
Kaushik Dayal.
\newblock Leading-order nonlocal kinetic energy in peridynamics for consistent
  energetics and wave dispersion.
\newblock {\em Journal of the Mechanics and Physics of Solids}, 105:235--253,
  2017.

\bibitem[DB06]{Kaushik-peri}
Kaushik Dayal and Kaushik Bhattacharya.
\newblock Kinetics of phase transformations in the peridynamic formulation of
  continuum mechanics.
\newblock {\em Journal of Mechanics and Physics of Solids}, 54:1811--1842,
  2006.

\bibitem[DBWW20]{dally2020cohesive}
Tim Dally, Carola Bilgen, Marek Werner, and Kerstin Weinberg.
\newblock Cohesive elements or phase-field fracture: Which method is better for
  dynamic fracture analyses?
\newblock In {\em Modeling and Simulation in Engineering-Selected Problems}.
  IntechOpen, 2020.

\bibitem[Eri75]{ericksen1975equilibrium}
Jerald~L Ericksen.
\newblock Equilibrium of bars.
\newblock {\em Journal of elasticity}, 5(3-4):191--201, 1975.

\bibitem[Esh80]{eshelby1980force}
JD~Eshelby.
\newblock The force on a disclination in a liquid crystal.
\newblock {\em Philosophical Magazine A}, 42(3):359--367, 1980.

\bibitem[FM06]{fuaciu2006longitudinal}
Cristian F{\u{a}}ciu and Alain Molinari.
\newblock On the longitudinal impact of two phase transforming bars. elastic
  versus a rate-type approach. part i: The elastic case.
\newblock {\em International journal of solids and structures},
  43(3-4):497--522, 2006.

\bibitem[FRMV17]{faye2017spherical}
Anshul Faye, Jos{\'e}~A Rodr{\'\i}guez-Mart{\'\i}nez, and KY~Volokh.
\newblock Spherical void expansion in rubber-like materials: The stabilizing
  effects of viscosity and inertia.
\newblock {\em International Journal of Non-Linear Mechanics}, 92:118--126,
  2017.

\bibitem[Gao97]{gao1997elastic}
Huajian Gao.
\newblock Elastic waves in a hyperelastic solid near its plane-strain
  equibiaxial cohesive limit.
\newblock {\em Philosophical magazine letters}, 76(5):307--314, 1997.

\bibitem[GLH{\etalchar{+}}19]{geelen2019phase}
Rudy~JM Geelen, Yingjie Liu, Tianchen Hu, Michael~R Tupek, and John~E Dolbow.
\newblock A phase-field formulation for dynamic cohesive fracture.
\newblock {\em Computer Methods in Applied Mechanics and Engineering},
  348:680--711, 2019.

\bibitem[HL85]{heidug1985thermodynamics}
W~Heidug and FK~Lehner.
\newblock Thermodynamics of coherent phase transformations in
  nonhydrostatically stressed solids.
\newblock {\em pure and applied geophysics}, 123(1):91--98, 1985.

\bibitem[KAILP17]{kumar2017some}
Aditya Kumar, Damian Aranda-Iglesias, and Oscar Lopez-Pamies.
\newblock Some remarks on the effects of inertia and viscous dissipation in the
  onset of cavitation in rubber.
\newblock {\em Journal of Elasticity}, 126(2):201--213, 2017.

\bibitem[KMB18]{kamensky2018hyperbolic}
David Kamensky, Georgios Moutsanidis, and Yuri Bazilevs.
\newblock Hyperbolic phase field modeling of brittle fracture: Part i—theory
  and simulations.
\newblock {\em Journal of the Mechanics and Physics of Solids}, 121:81--98,
  2018.

\bibitem[Kno02]{knowles2002impact}
James~K Knowles.
\newblock Impact-induced tensile waves in a rubberlike material.
\newblock {\em SIAM Journal on Applied Mathematics}, 62(4):1153--1175, 2002.

\bibitem[Lip14]{lipton2014dynamic}
Robert Lipton.
\newblock Dynamic brittle fracture as a small horizon limit of peridynamics.
\newblock {\em Journal of Elasticity}, 117(1):21--50, 2014.

\bibitem[LMTS{\etalchar{+}}18]{li2018variational}
Bin Li, Daniel Mill{\'a}n, Alejandro Torres-S{\'a}nchez, Benoit Roman, and
  Marino Arroyo.
\newblock A variational model of fracture for tearing brittle thin sheets.
\newblock {\em Journal of the Mechanics and Physics of Solids}, 119:334--348,
  2018.

\bibitem[LPM{\etalchar{+}}15]{li2015phase}
Bin Li, Christian Peco, Daniel Mill{\'a}n, Irene Arias, and Marino Arroyo.
\newblock Phase-field modeling and simulation of fracture in brittle materials
  with strongly anisotropic surface energy.
\newblock {\em International Journal for Numerical Methods in Engineering},
  102(3-4):711--727, 2015.

\bibitem[Mar06]{marder2006supersonic}
M~Marder.
\newblock Supersonic rupture of rubber.
\newblock {\em Journal of the Mechanics and Physics of Solids}, 54(3):491--532,
  2006.

\bibitem[MBB20]{mehrmashhadi2020validating}
Javad Mehrmashhadi, Mohammadreza Bahadori, and Florin Bobaru.
\newblock On validating peridynamic models and a phase-field model for dynamic
  brittle fracture in glass.
\newblock {\em Engineering Fracture Mechanics}, 240:107355, 2020.

\bibitem[PMB{\etalchar{+}}20]{peng20203d}
Xiaoyao Peng, Nithin Mathew, Irene~J Beyerlein, Kaushik Dayal, and Abigail
  Hunter.
\newblock A 3d phase field dislocation dynamics model for body-centered cubic
  crystals.
\newblock {\em Computational Materials Science}, 171:109217, 2020.

\bibitem[PZM{\etalchar{+}}20]{paul2020adaptive}
Karsten Paul, Christopher Zimmermann, Kranthi~K Mandadapu, Thomas~JR Hughes,
  Chad~M Landis, and Roger~A Sauer.
\newblock An adaptive space-time phase field formulation for dynamic fracture
  of brittle shells based on lr nurbs.
\newblock {\em Computational Mechanics}, 65(4):1039--1062, 2020.

\bibitem[RHA20]{rosakis2020inverse}
Phoebus Rosakis, Timothy~J Healey, and Ugur Alyanak.
\newblock The inverse-deformation approach to fracture.
\newblock {\em arXiv preprint arXiv:2006.16770}, 2020.

\bibitem[Ros95]{rosakis1995equal}
Phoebus Rosakis.
\newblock An equal area rule for dissipative kinetics of propagating strain
  discontinuities.
\newblock {\em SIAM Journal on Applied Mathematics}, 55(1):100--123, 1995.

\bibitem[SH20]{singh2020pseudomomentum}
H~Singh and JA~Hanna.
\newblock Pseudomomentum: origins and consequences.
\newblock {\em arXiv preprint arXiv:2007.06023}, 2020.

\bibitem[Sil00]{silling2000reformulation}
Stewart~A Silling.
\newblock Reformulation of elasticity theory for discontinuities and long-range
  forces.
\newblock {\em Journal of the Mechanics and Physics of Solids}, 48(1):175--209,
  2000.

\bibitem[TR14]{tupek2014extended}
MR~Tupek and R~Radovitzky.
\newblock An extended constitutive correspondence formulation of peridynamics
  based on nonlinear bond-strain measures.
\newblock {\em Journal of the Mechanics and Physics of Solids}, 65:82--92,
  2014.

\bibitem[Tru82]{truskinovskii1982equilibrium}
LM~Truskinovskii.
\newblock Equilibrium phase interfaces.
\newblock {\em Sov. Phys. Dokl.}, 27:551, 1982.

\bibitem[Tru93]{truskinovsky1993kinks}
Lev Truskinovsky.
\newblock Kinks versus shocks.
\newblock In {\em Shock induced transitions and phase structures in general
  media}, pages 185--229. Springer, 1993.

\bibitem[Tur97]{turteltaub1997viscosity}
Sergio Turteltaub.
\newblock Viscosity of strain gradient effects on the kinetics of propagating
  phase boundaries in solids.
\newblock {\em Journal of elasticity}, 46(1):53--90, 1997.

\bibitem[TV10]{trofimov2010shocks}
Evgeni Trofimov and Anna Vainchtein.
\newblock Shocks versus kinks in a discrete model of displacive phase
  transitions.
\newblock {\em Continuum Mechanics and Thermodynamics}, 22(5):317--344, 2010.

\bibitem[VTK17]{vidyasagar2017predicting}
Ananthan Vidyasagar, Wei~L Tan, and Dennis~M Kochmann.
\newblock Predicting the effective response of bulk polycrystalline
  ferroelectric ceramics via improved spectral phase field methods.
\newblock {\em Journal of the Mechanics and Physics of Solids}, 106:133--151,
  2017.

\bibitem[WA05]{weckner2005effect}
Olaf Weckner and Rohan Abeyaratne.
\newblock The effect of long-range forces on the dynamics of a bar.
\newblock {\em Journal of the Mechanics and Physics of Solids}, 53(3):705--728,
  2005.

\bibitem[YD10]{Lun_APL}
Lun Yang and Kaushik Dayal.
\newblock Formulation of phase-field energies for microstructure in complex
  crystal structures.
\newblock {\em Applied Physics Letters}, 96:081916, 2010.

\bibitem[ZKL16]{zhang2016variational}
Xiaoxuan Zhang, Andreas Krischok, and Christian Linder.
\newblock A variational framework to model diffusion induced large plastic
  deformation and phase field fracture during initial two-phase lithiation of
  silicon electrodes.
\newblock {\em Computer methods in applied mechanics and engineering},
  312:51--77, 2016.

\end{thebibliography}

\end{document}